\documentclass[11pt]{article}

\usepackage[a4paper]{geometry}
\usepackage[show]{ed}   

\usepackage{amsthm}   
\usepackage{xcolor}
\usepackage{url}
\usepackage{tabu}
\usepackage{bm}

\newtheorem{theorem}{Theorem}[section]
\newtheorem{assumption}[theorem]{Assumption}
\newtheorem{corollary}[theorem]{Corollary}
\newtheorem{definition}[theorem]{Definition}

\newtheorem{lemma}[theorem]{Lemma}

\newtheorem{proposition}[theorem]{Proposition}
\newtheorem{remark}[theorem]{Remark}

\usepackage{}
\usepackage{graphicx}
\usepackage{multirow,bigstrut}
\usepackage{amsmath,amsfonts,amssymb}
\usepackage{mathrsfs}
\usepackage{color}
\usepackage{upgreek}
\usepackage{bm}
\usepackage{verbatim}
\usepackage{mathtools}
\usepackage{calligra}
\usepackage[T1]{fontenc}
\usepackage{subfigure}
\usepackage{latexsym}
\usepackage{color}
\usepackage{soul}
\usepackage{authblk}

\usepackage{verbatim}
\usepackage{exscale}
\usepackage{relsize}

\usepackage{tikz}
\usetikzlibrary{shapes,arrows}

\usepackage{authblk}
\usepackage{hyperref}

\numberwithin{equation}{section}

\newcommand{\cB}{\color{black}}

\newcommand{\cR}{\color{black}}

\usepackage{}
\usepackage{graphicx}
\usepackage{multirow,bigstrut}
\usepackage{amsmath,amsfonts,amssymb}
\usepackage{mathrsfs}
\usepackage{color}
\usepackage{upgreek}
\usepackage{bm}
\usepackage{verbatim}
\usepackage{mathtools}
\usepackage{calligra}
\usepackage[T1]{fontenc}
\usepackage{subfigure}
\usepackage{pgfplots}
\usepackage{amsmath,accents}
\usepackage{cleveref}

\usepackage{diagbox}
\usepackage{booktabs}
\usepackage{multicol}
\usepackage{subfloat}
\usepackage{tikz}

\usepackage{tikz}
  \usetikzlibrary{calc,shapes,decorations,decorations.text, mindmap,shadings,patterns,matrix,arrows,intersections,automata,backgrounds}
  
\usetikzlibrary{shapes,arrows}

\usepackage{verbatim}
\usepackage{exscale}
\usepackage{relsize}

\DeclareMathAlphabet{\mathdutchcal}{U}{dutchcal}{m}{n}
\SetMathAlphabet{\mathdutchcal}{bold}{U}{dutchcal}{b}{n}
\DeclareMathAlphabet{\mathdutchbcal}{U}{dutchcal}{b}{n}

\usepackage[T1]{fontenc}
\DeclareFontFamily{T1}{calligra}{}
\DeclareFontShape{T1}{calligra}{m}{n}{<->s*[1.44]callig20}{}
\DeclareMathAlphabet\mathcalligra   {T1}{calligra} {m} {n}
\DeclareMathAlphabet\mathzapf       {T1}{pzc} {mb} {it}
\DeclareMathAlphabet\mathchorus     {T1}{qzc} {m} {n}
\DeclareMathAlphabet\mathrsfso      {U}{rsfso}{m}{n}

\makeatletter
\tikzset{
        hatch distance/.store in=\hatchdistance,
        hatch distance=5pt,
        hatch thickness/.store in=\hatchthickness,
        hatch thickness=5pt
        }
\pgfdeclarepatternformonly[\hatchdistance,\hatchthickness]{north east hatch}
    {\pgfqpoint{-1pt}{-1pt}}
    {\pgfqpoint{\hatchdistance}{\hatchdistance}}
    {\pgfpoint{\hatchdistance-1pt}{\hatchdistance-1pt}}%
    {
        \pgfsetcolor{\tikz@pattern@color}
        \pgfsetlinewidth{\hatchthickness}
        \pgfpathmoveto{\pgfqpoint{0pt}{0pt}}
        \pgfpathlineto{\pgfqpoint{\hatchdistance}{\hatchdistance}}
        \pgfusepath{stroke}
    }
\makeatother

\newcommand{\vertiii}[1]{{\left\vert\kern-0.25ex\left\vert\kern-0.25ex\left\vert #1 
    \right\vert\kern-0.25ex\right\vert\kern-0.25ex\right\vert}}

\newcommand{\cm}[1]{{\color{black}{#1}}}
\newcommand{\chupeng}[1]{{\color{black}{#1}}}

\newcommand{\cpm}[1]{{\color{black}{#1}}}

\newcommand{\alber}[1]{\textcolor{black}{#1}}

\newcommand{\rs}[1]{\textcolor{black}{#1}}

\newcommand{\RS}[1]{\textcolor{black}{#1}}

\begin{document}


\title{Two-level Restricted Additive Schwarz Preconditioner based on Multiscale Spectral Generalized FEM for Heterogeneous Helmholtz Problems}


\date{} 

\author[1]{Chupeng Ma}
\author[2]{Christian Alber}
\author[2]{Robert Scheichl}
\author[3]{Yongwei Zhang}
\affil[1]{School of Sciences, Great Bay University, 523000 Dongguan, China}
\affil[2]{Institute for Mathematics and Interdisciplinary Center for Scientific Computing, Heidelberg University, Im Neuenheimer Feld 205, Heidelberg 69120, Germany}
\affil[3]{School of Mathematics and Statistics, Zhengzhou University, Zhengzhou University, 450001 Zhengzhou, China}
\maketitle

\begin{abstract}
We present and analyze a two-level restricted additive Schwarz (RAS) preconditioner for heterogeneous Helmholtz problems, based on a multiscale spectral generalized finite element method (MS-GFEM) proposed in [C. Ma, C. Alber, and R. Scheichl, SIAM. J. Numer. Anal., 61 (2023), pp. 1546--1584]. The preconditioner uses local solves with impedance boundary conditions, and a global coarse solve based on the MS-GFEM approximation space constructed from local eigenproblems. It is derived by first formulating MS-GFEM as a Richardson iterative method, and without using an oversampling technique, reduces to the preconditioner recently proposed and analyzed in [Q. Hu and Z.Li, arXiv 2402.06905]. 

We prove that both the Richardson iterative method and the preconditioner used within GMRES converge at a rate of $\Lambda$ under some reasonable conditions, where $\Lambda$ denotes the error of the underlying MS-GFEM \rs{approximation}. Notably, the convergence proof of GMRES does not rely on the ‘Elman theory'. An exponential convergence property of MS-GFEM, resulting from oversampling, ensures that only a few iterations are needed for convergence with a small coarse space. Moreover, the convergence rate $\Lambda$ is not only independent of the fine-mesh size $h$ and the number of subdomains, but decays with increasing wavenumber $k$. In particular, in the constant-coefficient case, with $h\sim k^{-1-\gamma}$ for some $\gamma\in (0,1]$, it holds that $\Lambda \sim k^{-1+\frac{\gamma}{2}}$. We present extensive numerical experiments to illustrate the performance of the preconditioner, including 2D and 3D benchmark geophysics tests, and a high-contrast coefficient example arising in applications. 
\end{abstract}

\noindent
{\bf MSC Classes:} {\tt 65N30, 65N55, 65F10}
\\[2ex]

\section{Introduction}
The Helmholtz equation is a basic mathematical model for describing wave phenomena in acoustics and electromagnetics. The numerical solution of the Helmholtz equation at high frequency is notoriously difficult, for two main reasons. First, due to the so-called \textit{pollution effect}, standard discretization methods typically \alber{need} a more restrictive condition on the mesh-size $h$ than $h\sim k^{-1}$ ($k$ denoting the wavenumber) \cite{melenk2011wavenumber,du2015preasymptotic}, namely, a higher mesh resolution than needed for a reasonable representation of the solution. With $k\rightarrow \infty$, this leads to a huge number of unknowns in the linear system to be solved, making it not amenable to direct solvers. Secondly, the resulting linear system is non-Hermitian and highly indefinite at high frequency, for which classical iterative solvers and preconditioning techniques exhibit slow convergence \cite{ernst2012difficult}. Apart from the high-frequency challenge, practical applications such as seismic full-waveform inversion (FWI) \cite{virieux2017introduction} often involve \rs{the solution of} Helmholtz equations in strongly heterogeneous media, possibly for a large number of sources. The interaction of high-frequency waves with small-scale heterogeneities renders the numerical solution of heterogeneous Helmholtz problems even more challenging. Moreover, multiple-source simulations call for fast solution of a single source.

Considerable research has been performed on iterative solvers for the Helmholtz equation; see literature surveys in \cite{graham2020domain,gander2019class,zepeda2016method}. Two successful preconditioners in this regard are sweeping preconditioners \cite{engquist2011sweeping,engquist2011sweeping-2,zepeda2016method,taus2020sweeps} and the shifted Laplace preconditioner \cite{erlangga2006novel,van2007spectral}. In this paper, we are concerned with parallel iterative solvers based on domain decomposition methods (DDMs). While one-level DDMs, based on the parallel solution of a set of local problems at each iteration, are relatively simple and have demonstrated \rs{significant} success \cite{dolean2015introduction,gander2002optimized,gander2007optimized,gong2021domain,gong2022convergence,gong2023convergence,graham2020domain}, we focus on two-level methods with an additional coarse-space component. The coarse component is essential to achieve scalability with respect to the number of subdomains, and robustness with respect to problem parameters. Compared with Poisson-type equations, the main challenge of two-level DDMs for Helmholtz equations is that to achieve a $k$-independent convergence, the coarse space should be reasonably `fine' to resolve the wave nature of the solution. Using the one with superior approximation properties is an effective strategy to reduce the coarse space size.

Earlier works on two-level DDMs for the Helmholtz equation include \cite{cai1992domain}, \chupeng{which uses} the classical coarse space based on piecewise polynomials on a coarse grid, and \cite{kimn2007restricted,farhat2000two}, \chupeng{which uses} coarse spaces based on plane waves. We refer to \cite{graham2017domain,sieburgh2024coarse} for other related works in this direction. In recent years, spectral coarse spaces based on local generalized eigenproblems have been popular in designing coefficient-robust DD preconditioners for multiscale PDEs; see, e.g., \cite{galvis2010domain,galvis2010domain-2,heinlein2019adaptive,spillane2014abstract,spillane2023toward,nataf2024coarse}. These coarse spaces are built by selecting a certain number of eigenfunctions per subdomain based on a user-specific tolerance for the eigenvalues. A particular spectral coarse space relevant to our work is the GenEO (Generalized Eigenproblems in the Overlap) coarse space \cite{spillane2014abstract}. The concept of spectral coarse spaces was later extended to the Helmholtz equation to construct wavenumber-robust preconditioners. Two representative coarse spaces in this regard are the DtN and GenEO type (\rs{termed} ‘H-GenEO') coarse spaces \cite{bootland2021comparison,conen2014coarse}, which are based on local eigenproblems suitably adapted to the Helmholtz equation. While the two coarse spaces show promising numerical results, a rigorous convergence theory remains out of reach for them. In \cite{bootland2023overlapping}, with another GenEO-type coarse space defined with positive definite operators, the convergence of a two-level Schwarz preconditioner within GMRES was proved for Helmholtz-type equations. The convergence results have been improved in \cite{dolean2024improvements} recently.


More recently, there has been a growing interest in developing multiscale/spectral coarse spaces based two-level Schwarz preconditioners for the Helmholtz equation \cite{dolean2024schwarz,hu2024novel,fu2024edge,lu2024two}. These coarse spaces include GenEO-type coarse spaces \cite{dolean2024improvements,hu2024novel}, and multiscale approximation spaces \cite{lu2024two,fu2024edge} based on the localized orthogonal decomposition (LOD) technique \cite{maalqvist2014localization} and an edge multiscale method \cite{fu2021wavelet}. In all these works, it was proved that the GMRES or Richardson iterative solvers applied with the proposed preconditioners converge uniformly independent of the wavenumber (and other important parameters) under certain assumptions on the subdomain and coarse-space sizes. In particular, we note that \cite{hu2024novel} \chupeng{use} a GenEO-type coarse space with local eigenproblems posed on discrete Helmholtz-harmonic spaces, corresponding to the multiscale approximation space proposed in \cite{ma2023wavenumber} (except without oversampling). The resulting two-level restricted additive Schwarz (RAS) preconditioner was proved to converge uniformly, provided that the subdomain size $H$ and the eigenvalue tolerance $\rho$ satisfy $H\lesssim k^{-1}$ and $\rho \lesssim k^{-1}$. Also, we note that a similar coarse space was recently proposed in \cite{nataf2024coarse} for general domain decomposition methods in a fully algebraic setting. While writing this paper, we found two new works on two-level Schwarz preconditioning for Helmholtz equations \cite{galkowski2025convergence,graham2025two}. The authors in \cite{galkowski2024convergence} proved that under certain assumptions on approximation properties of a coarse space, which are satisfied by the coarse spaces in \cite{hu2024novel,fu2024edge,lu2024two}, the preconditioned matrix of a particular two-level Schwarz method is close to the identity. In \cite{graham2025two}, the authors studied two-level Schwarz preconditioners with piecewise-polynomial coarse spaces, and identified precise conditions for the preconditioned GMRES to converge with iteration number that grows at most like $(\log k)^{2}$.  


While there have been significant developments in the theory of two-level Schwarz methods based on multiscale/spectral coarse spaces for the Helmholtz equation, several issues remain to be addressed. First, although multiscale or spectral coarse spaces were originally designed for PDEs with strongly heterogeneous coefficients, most of recent works (e.g., \cite{hu2024novel,fu2024edge,lu2024two}) focus only on the constant-coefficient Helmholtz equation, and it is unclear whether the theories developed apply to the case of heterogeneous coefficients. Second, regarding spectral coarse spaces, existing convergence theories require the eigenvalue tolerance to satisfy similar conditions as $\rho \lesssim k^{-1}$ \cite{bootland2023overlapping,dolean2024improvements,hu2024novel}. Such conditions become more restrictive as the wavenumber increases. An important question naturally arises: how fast do the eigenvalues decay? The faster the eigenvalues decay, the less eigenfunctions are needed to satisfy the criterion $\rho \lesssim k^{-1}$, and thus the smaller the coarse problem is. However, this question was rarely investigated theoretically in related works. Finally, a ‘standard' tool for analysing the convergence of GMRES is the ‘Elman theory' \cite{eisenstat1983variational}, on which most existing convergence analyses of DDMs for the Helmholtz equation \chupeng{rely} \cite{bootland2023overlapping,cai1992domain,dolean2024schwarz,gong2021domain,hu2024novel,lu2024two}. However, the bounds \chupeng{on} the convergence rates of GMRES are far from sharp (see Chapter 11 of \cite{toselli2006domain}). Indeed, given the high cost of multiscale/spectral coarse spaces, it is reasonable to expect that a coarse space (if incorporated suitably) may be able to do more in terms of convergence speed than only achieving scalability and robustness. However, even if such a coarse space exists, the ‘Elman theory', in general, is unable to predict this. This paper aims to address \RS{all} these issues.


In this paper, we present and analyze a two-level hybrid RAS preconditioner for Helmholtz problems with strongly heterogeneous coefficients, based on the Multiscale Spectral Generalized Finite Element Method (MS-GFEM) \cite{ma2023wavenumber}. As a multiscale model reduction method, MS-GFEM builds a coarse model (coarse space) from selected eigenfunctions of GenEO-type local eigenproblems posed on Helmholtz-harmonic spaces. By solving these eigenproblems in oversampled subdomains, MS-GFEM achieves local exponential convergence (with respect to the number of local spectral basis functions). We derive the preconditioner by first formulating MS-GFEM as a Richardson iterative method following \cite{strehlow2024fast}, where symmetric positive definite problems were considered. However, we use impedance conditions instead of Dirichlet conditions as in \cite{strehlow2024fast} to construct local problems, leading to a two-level ORAS (‘Optimised Restricted Additive Schwarz’) type preconditioner. Except for oversampling, this preconditioner exactly corresponds to the one proposed in \cite{hu2024novel}. Following \cite{strehlow2024fast}, we rigorously prove that both the Richardson iterative method and the preconditioner used within GMRES converge at a rate of $\Lambda$, where $\Lambda$ denotes the error of the underlying MS-GFEM, provided that similar resolution conditions as in \cite{hu2024novel} -- required for quasi-optimality of the underlying MS-GFEM -- are satisfied. Thanks to an exponential decay of the eigenvalues and a resulting exponential convergence of MS-GFEM established in \cite{ma2023wavenumber}, our result ensures that only a few iterations are needed for convergence with a small coarse space. Notably, we prove the convergence of GMRES without using the ‘Elman theory', and show that the convergence rate $\Lambda$ decays with increasing wavenumber $k$ as $k^{-\alpha}$ for some $\alpha > 0$. In particular, in the constant-coefficient case, with $h\sim k^{-1-\gamma}$ for some $0<\gamma\leq 1$, we show that 
\begin{equation*}
\Lambda \sim k^{-1+\frac{\gamma}{2}}, \quad \text{or equivalently},\quad \Lambda \sim \rho^{1-\frac{\gamma}{2}},   
\end{equation*}
where $\rho$ denotes the eigenvalue tolerance. In addition, we prove that a variant of the preconditioner -- with Dirichlet conditions for the local solves -- has a better convergence rate $\Lambda\sim k^{-1}$ under the same conditions, which is therefore independent of $\gamma$.

While the presented preconditioner has already been proposed and analysed in \cite{hu2024novel}, there are several differences between our work and \cite{hu2024novel}, especially in terms of methodology and convergence results. First, in contrast to \cite{hu2024novel}, we consider heterogeneous Helmholtz problems, and our results hold for the general case of $L^{\infty}$-coefficients. Secondly, the preconditioner in this paper is derived naturally from MS-GFEM, whereas the one in \cite{hu2024novel} is defined directly by adding a separate coarse space to \rs{a} one-level preconditioner (similarly as in \cite{fu2024edge,lu2024two}). As a result, our preconditioner exhibits a more direct and clear connection with the underlying multiscale method, which plays an important role in the convergence analysis. Thirdly, we prove the convergence of GMRES based on the convergence of the Richardson iterative method and a residual minimization property of GMRES, instead of using the ‘Elman theory' as in \cite{hu2024novel}. This yields a much sharper bound with a precise characterization for the convergence rate of GMRES. Finally, we use oversampling in the definition of the preconditioner, leading to a provably exponential decay of the eigenvalues. The oversampling technique significantly reduces the coarse space size for a fixed eigenvalue tolerance.  

The rest of this paper is organized as follows. \Cref{sec:MS-GFEM} gives a detailed description of the MS-GFEM for solving discretized Helmholtz problems. The two-level Schwarz preconditioner with its convergence theory is presented in \cref{sec:preconditioner} . Numerical experiments are given in \cref{sec:numerical results}. Finally, we give a few concluding remarks in \cref{sec:conclusion}.

\section{MS-GFEM for discretized Helmholtz problems}\label{sec:MS-GFEM}
In this section, we present MS-GFEM for solving discretized Helmholtz problems following \cite{ma2023wavenumber}, with a focus on its convergence theory. 

\subsection{Problem formulation and discretization}
Let $\Omega\subset \mathbb{R}^{d}$ $(d=2,3)$ be a bounded Lipschitz domain with a polygonal boundary $\partial \Omega$. We consider the following heterogeneous Helmholtz equation with wavenumber $k>0$: Find $\cm{u^{\mathdutchcal{e}}}:\Omega\rightarrow\mathbb{C}$ such that
\begin{equation}\label{eq:1-1}
\left\{
\begin{array}{lll}
{\displaystyle -{\rm div}(A\nabla \cm{u^{\mathdutchcal{e}}}) - k^{2}V^{2}\cm{u^{\mathdutchcal{e}}}= f\,\quad {\rm in}\;\, \Omega, }\\[2mm]
{\displaystyle \qquad \,\;A\nabla \cm{u^{\mathdutchcal{e}}}\cdot {\bm n} - {\rm i}k\beta \cm{u^{\mathdutchcal{e}}}=g  \quad \,{\rm on}\;\,\partial \Omega,}
\end{array}
\right.
\end{equation}
where ${\bm n}$ denotes the outward unit normal to $\partial \Omega$. For simplicity, we only consider the pure impedance boundary condition in this paper. Extension of the results below to the case of mixed boundary conditions is straightforward; see \cite{ma2023wavenumber}. Throughout this paper, we make the following assumptions on the data $A$, $V$, $\beta$, $f$, and $g$:  

\begin{assumption}\label{ass:problem-setting} 
\begin{itemize}
\item[(i)] $A \in L^{\infty}(\Omega,\mathbb{R}_{\rm sym}^{d\times d})$ is uniformly elliptic: there exists $0< a_{\rm min} < a_{\rm max}<\infty$ such that
\begin{equation}\label{eq:1-1-0}
a_{\rm min} |{\bm \xi}|^{2} \leq A({\bm x}){\bm \xi}\cdot{\bm \xi} \leq a_{\rm max}  |{\bm \xi}|^{2},\quad \forall {\bm \xi}\in \mathbb{R}^{d},\quad {\bm x} \in\Omega;
\end{equation}
\item[(ii)] $V \in L^{\infty}(\Omega)$ and there exists $0<V_{\rm min}<V_{\rm max}<\infty$ such that $V_{\rm min}\leq V({\bm x})\leq V_{\rm max}$ for a.e. ${\bm x}\in \Omega$;

\vspace{1ex}
\item[(iii)] $f\in L^{2}(\Omega)$, $g\in L^{2}(\partial \Omega)$ and $\beta\in L^{\infty}(\partial \Omega)$. \cm{Moreover, $\beta({\bm x})>0$ \cpm{for a.e.} ${\bm x}\in \partial \Omega$}.
\end{itemize}
\end{assumption}

Let $\mathcal{B}: H^{1}(\Omega)\times H^{1}(\Omega) \rightarrow \mathbb{C}$ be the sesquilinear form defined by
\begin{equation}
\mathcal{B}(u,v) = \int_{\Omega} \big(A\nabla u\cdot {\nabla \overline{v}}-k^{2}V^{2}u\overline{v}\big) \,d{\bm x} - {\rm i}k\int_{\partial \Omega} \beta u\overline{v}\,d{\bm s}\quad \forall u,\,v\in H^{1}(\Omega).
\end{equation}
The weak formulation of problem \eqref{eq:1-1} is as follows: Find $u^{\mathdutchbcal{\mathdutchcal{e}}}\in H^{1}(\Omega)$ such that 
\begin{equation}\label{continous-HelmholtzPDE}
\mathcal{B}(u^{{\mathdutchcal{e}}},v) = F(v):=\int_{\partial \Omega} g\overline{v} \,d{\bm s}+ \int_{\Omega}f\overline{v} \,d{\bm x}\quad \forall v\in H^{1}(\Omega).
\end{equation}

We now consider the finite element approximation for problem \eqref{continous-HelmholtzPDE}. Let $\{\tau_{h}\}_{h}$ be a family of shape-regular triangulations of $\Omega$ with $h$ denoting the mesh size. Let $U_{h}(\Omega)\subset H^{1}(\Omega)$ be a Lagrange finite element space on $\tau_h$ consisting of continuous piecewise polynomials of degree $p$. The standard finite element discretization of problem \eqref{continous-HelmholtzPDE} is to find $u^{\mathdutchcal{e}}_{h}\in U_{h}(\Omega)$ such that
\begin{equation}\label{fineFE_problem}
\mathcal{B}(u^{\mathdutchcal{e}}_{h},v_{h}) = F(v_{h})\qquad \forall v_{h}\in U_{h}(\Omega).
\end{equation}
Denote by $m$ the dimension of $U_{h}(\Omega)$ and let $\{\phi_{k}\}_{k=1}^{m}$ be a basis for $U_{h}(\Omega)$. Then we can write \eqref{fineFE_problem} as the linear system:
\begin{equation}\label{linear_system}
 {\bf B}{\bf u} = {\bf F},   
\end{equation}
where \chupeng{${\bf B} = (b_{lj})\in \mathbb{C}^{m\times m}$ and ${\bf F} = (F_j)\in \mathbb{C}^{m}$, with $b_{lj} = \mathcal{B}(\phi_l,\phi_j)$ and $F_j = F(\phi_j)$}. 

For later use, let us introduce some local sesquilinear forms. For any subdomain $\omega\subset \Omega$ and any $u$, $v\in H^{1}(\omega)$, we define 
\begin{equation}\label{localforms}
\begin{array}{cc}
{\displaystyle \mathcal{B}_{\omega}(u,v) = \int_{\omega} \big(A\nabla u\cdot {\nabla \overline{v}}-k^{2}V^{2}u\overline{v}\big) \,d{\bm x} - {\rm i}k\int_{ \partial \omega\cap \partial \Omega} \beta u\overline{v}\,d{\bm s},}\\[3mm]
{\displaystyle \mathcal{B}_{\omega,{\rm imp}}(u,v) = \mathcal{B}_{\omega}(u,v) - {\rm i}k\int_{ \partial \omega\cap \Omega} V u\overline{v}\,d{\bm s},}\\[3mm]
{\displaystyle \mathcal{A}_{\omega,k}(u,v) = \int_{\omega} \big(A\nabla u\cdot {\nabla \overline{v}}+k^{2}V^{2}u\overline{v}\big)\,d{\bm x}, }\\[3mm]
{\displaystyle \mathcal{A}_{\omega}(u,v) = \int_{\omega} A\nabla u\cdot {\nabla \overline{v}}\,d{\bm x},}
\end{array}
\end{equation}
and set
\begin{equation}\label{eq:1-4}
\begin{array}{lll}
{\displaystyle \Vert u \Vert_{\mathcal{A},\omega} =  \sqrt{\mathcal{A}_{\omega}(u,u)},\quad \Vert u \Vert_{\mathcal{A},\omega,k} = \sqrt{\mathcal{A}_{\omega,k}(u,u)}.}
\end{array}
\end{equation}
When $\omega=\Omega$, we simply write $\Vert u \Vert_{\mathcal{A}}$ and $\Vert u \Vert_{\mathcal{A},k}$ for the norms. The sesquilinear form $\mathcal{B}(\cdot,\cdot)$ has the following boundedness property under $\Vert\cdot\Vert_{\mathcal{A},k}$ (see \cite{melenk1995generalized}): there exists a $k$-independent constant $C_{\mathcal{B}}>0$ such that
\begin{equation}\label{boundedness-estimate}
|\mathcal{B}(u,v)|\leq C_{\mathcal{B}} \Vert u \Vert_{\mathcal{A},k} \Vert v \Vert_{\mathcal{A},k},\quad \forall u,\,v\in H^{1}(\Omega).
\end{equation}
Throughout this paper, we use the notation $a\lesssim b$ ($a\gtrsim b$) to denote the inequality $a\leq Cb $ (resp. $a\geq Cb$), where $C$ is a positive constant independent of all key parameters in our analysis, such as the mesh-size $h$, the wavenumber $k$, and the size of the oversampling domains below, but may depend on the coefficient bounds in Assumption~\ref{ass:problem-setting}. We write $a\sim b$ if $a\lesssim b$ and $b\gtrsim a$.

We make the following assumptions concerning the solvability of problem \eqref{continous-HelmholtzPDE} and its discretization \eqref{fineFE_problem}.
\begin{assumption}\label{ass:continuous-solvablity}
For each $f\in L^{2}(\Omega)$ and $g\in L^{2}(\partial \Omega)$, the problem \eqref{continous-HelmholtzPDE} has a unique solution $u^{\mathdutchbcal{\mathdutchcal{e}}}\in H^{1}(\Omega)$. Moreover, there exists $C_{\mathtt{stab}}(k)>0$ depending polynomially on $k$, such that the solution $u^{\mathdutchbcal{\mathdutchcal{e}}}$ of problem \eqref{continous-HelmholtzPDE} with $g=0$ satisfies
\begin{equation}
\Vert u^{\mathdutchbcal{\mathdutchcal{e}}} \Vert_{\mathcal{A},k}\leq  C_{\mathtt{stab}}(k) \Vert f\Vert_{L^{2}(\Omega)}.   
\end{equation}
\end{assumption}
\begin{remark}\label{rem:well-posedness}
\rs{In two dimensions, the} well-posedness of problem \eqref{continous-HelmholtzPDE} was proved under Assumption~\ref{ass:problem-setting}, whereas in \rs{three dimensions}, it requires an additional assumption that $A$ is piecewise Lipschitz; see \cite[Lemma 2.4]{gong2021domain}. Concerning the stability constant $C_{\mathtt{stab}}(k)$, we assume that it depends polynomially on $k$ to exclude some extreme ‘trapping' cases \cite{chandler2020high} where $C_{\mathtt{stab}}(k)$ depends exponentially on $k$. In some nontrapping cases, \cite[Theorem 3.10]{gong2021domain} \rs{shows} that $C_{\mathtt{stab}}(k)$ is independent of $k$. In particular, if $A = I$, $V=1$, $\beta=1$, and $\Omega$ is a star-shaped domain with respect to a ball of radius $\sim L$, it was proved that $C_{\mathtt{stab}}(k) \sim L$\rs{; see} \cite[Theorem 2.7]{graham2020domain}.
\end{remark}

\begin{assumption}\label{ass:global_stability}
There exists $h_{0}>0$ such that for any $0<h<h_{0}$, the estimate 
\begin{equation}\label{global_discrete_infsup} 
\inf_{v_{h}\in U_{h}(\Omega)} \sup_{0\neq u_{h}\in U_{h}(\Omega)} \frac{|\mathcal{B}(v_{h}, u_{h})|}{ \Vert v_{h}\Vert_{\mathcal{A},k} \Vert u_{h}\Vert_{\mathcal{A},k}}  \gtrsim \frac{1}{1 + kC_{\mathtt{stab}}(k)}, 
\end{equation}
holds with $C_{\mathtt{stab}}(k)$ defined in Assumption~\ref{ass:continuous-solvablity}.
\end{assumption}
\begin{remark}\label{rem:stab_constant}
In the constant-coefficient case, it was proved that the estimate \eqref{global_discrete_infsup} holds if $k^{p+1}h^{p}$ is sufficiently small; see \cite{graham2020domain}. However, for the well-posedness of problem \eqref{fineFE_problem}, one only needs the weaker condition that $k^{2p+1}h^{2p}$ is sufficiently small; see \cite{du2015preasymptotic,wu2014pre}. In the heterogeneous-coefficient case, estimate \eqref{global_discrete_infsup} was proved in \cite[Lemma 3.16]{gong2021domain} under a condition on the so-called ‘adjoint approximability'; see also \cite[Theorem 3.2]{melenk2011wavenumber}. \rs{In} \cite[Theorem 2]{schatz1996some} \rs{it was also shown} that such a condition is valid for general heterogeneous coefficients if the mesh-size is sufficiently small, and a $k$-explicit characterization of the threshold mesh-size $h_0$ was given in \cite{chaumont2020wavenumber,graham2020stability} for some special problem settings.  
\end{remark}


\subsection{MS-GFEM}
We construct the MS-GFEM within the framework of overlapping DDMs. To this end, we first decompose the domain $\Omega$ into a set of overlapping subdomains $\{\omega_{i}\}_{i=1}^{N}$ and assume that each $\omega_i$ is a union of mesh elements of $\tau_h$. For each $\omega_i$, we introduce an oversampling domain $\omega_{i}^{\ast} \supset \omega_i$ formed by adding several layers of adjacent mesh elements to $\omega_i$, which we assume is a Lipschitz domain. We define the following local finite element spaces:
\begin{equation}
\begin{array}{lll}
{\displaystyle  {U}_{h}(\omega^{\ast}_{i}) = \big\{v_{h}|_{\omega^{\ast}_{i}}\;:\; v_{h}\in U_{h}\big\},}\\[2mm]
{\displaystyle U_{h,0}(\omega^{\ast}_i)= \big\{v_{h}\in {U}_{h}(\omega^{\ast}_{i}):\;v_{h} = 0 \;\;{\rm on}\;\, \partial \omega^{\ast}_{i} \cap \Omega \big\}, }\\[2mm]
{\displaystyle U_{h,\mathcal{B}}(\omega^{\ast}_i)= \big\{u_{h}\in U_{h}(\omega^{\ast}_i)\;:\; \mathcal{B}_{\omega^{\ast}_{i}}(u_{h},v_{h}) = 0\;\, \forall v_{h}\in  U_{h,0}(\omega^{\ast}_i)\big\}.}
\end{array}
\end{equation}
$U_{h,\mathcal{B}}(\omega^{\ast}_i)$ is a discretized Helmholtz-harmonic space. 


Next we define the local particular functions for MS-GFEM. Denote by $\mathcal{N}^{h} = \{{\bf x}_q: q=1,\cdots, m\}$ the \rs{set of nodes} associated with the basis $\{\phi_i\}_{i=1}^{m}$ of $U_{h}(\Omega)$. For each $1\leq i\leq N$, we introduce a node-wise (zero) extension operator $E_{i}:U_{h}(\omega_i^{\ast})\rightarrow U_{h}(\Omega)$ defined for \rs{each} ${\bf x}_q \in \mathcal{N}^{h}$ \rs{as follows:}
\begin{equation}
(E_{i}v_{h})({\bf x}_q) = 
\left\{
    \begin{array}{cc}
         v_{h}({\bf x}_q)\quad {\rm if}\;\;{\bf x}_q\in \overline{\omega_i^{\ast}},  \\[2mm]
         0  \quad \quad \quad {\rm otherwise}.
    \end{array}
    \right.
\end{equation}
On each oversampling domain $\omega_i^{\ast}$, we consider the following local problem: Find $\psi_{h,i}\in {U}_{h}(\omega^{\ast}_{i})$ such that
\begin{equation}\label{localHelm_problem}
\mathcal{B}_{\omega_{i}^{\ast},{\rm imp}}(\psi_{h,i},v_{h}) = F(E_{i}v_{h})\quad \forall v_{h}\in {U}_{h}(\omega^{\ast}_{i}),
\end{equation}
where the local sesquilinear form $B_{\omega_i^{\ast},{\rm imp}}(\cdot,\cdot)$ is defined in \eqref{localforms}. Note that \rs{this imposes a} homogeneous impedance boundary condition 
\begin{equation}\label{local_IMBC}
A\nabla \psi_{h,i}\cdot {\bm n} - {\rm i}kV \psi_{h,i} = 0 \quad \text{on}\;\; \partial\omega_i^{\ast}\cap \Omega    
\end{equation}
for $\psi_{h,i}$ \rs{on the artificial boundary of} the variational problem \eqref{localHelm_problem}. Combining \eqref{fineFE_problem} and \eqref{localHelm_problem}, one can easily see that $u^{e}_{h}|_{\omega_i^{\ast}} - \psi_{h,i}\in U_{h,\mathcal{B}}(\omega^{\ast}_i)$, i.e., $u^{e}_{h}|_{\omega_i^{\ast}} - \psi_{h,i}$ is discrete Helmholtz-harmonic.
\begin{remark}
In \cite{ma2023wavenumber}, the local particular function was defined in a slightly different way: Find $\psi_{h,i}\in {U}_{h}(\omega^{\ast}_{i})$ such that 
\begin{equation}\label{localHelm_Original}
\mathcal{B}_{\omega_{i}^{\ast},{\rm imp}}(\psi_{h,i},v_{h}) = F(\widetilde{v_{h}})\quad \forall v_{h}\in {U}_{h}(\omega^{\ast}_{i}),
\end{equation}  
where $\widetilde{v_{h}}$ denotes the zero extension of $v_{h}$ outside $\omega_i^{\ast}$. Here, in order to formulate MS-GFEM as an iterative solver, we define it in the form of \eqref{localHelm_problem} following \cite{gong2023convergence,hu2024novel}. 
\end{remark}
\begin{remark}\label{rem:DBC_localsolution}
Alternatively, we can define the local particular function $\psi_{h,i} \in {U}_{h,0}(\omega^{\ast}_{i})$ such that 
\begin{equation}\label{local_DBC}
\mathcal{B}_{\omega_{i}^{\ast},{\rm imp}}(\psi_{h,i},v_{h}) = F(v_{h})\quad \forall v_{h}\in {U}_{h,0}(\omega^{\ast}_{i}),    
\end{equation}
i.e., we use a zero Dirichlet boundary condition for $\psi_{h,i}$. The local problem \eqref{local_DBC} is well defined if ${\rm diam}(\omega^{\ast}_i) \lesssim k^{-1}$. The benefit of using this definition of $\psi_{h,i}$ is that we can get rid of the term $(kh)^{-1/2}$ in the stability estimate (Lemma~\ref{lem:error_local_Helmholtz}) and in the resulting final error estimate (Corollary~\ref{cor:final_approx_result}) below; see Remark~\ref{rem:DBC_local_stability} and Remark~\ref{rem:DBC_final_error_result}.
\end{remark}

We proceed to build the local approximation spaces for MS-GFEM. To do this, we first introduce a partition of unity $\{ \chi_{i} \}_{i=1}^{N}$ subordinate to the open cover $\{\omega_{i}\}_{i=1}^{N}$ with the properties
\begin{equation}\label{eq:2-0-1}
\begin{array}{lll}
{\displaystyle 0\leq \chi_{i}\leq 1,\quad {\rm supp}(\chi_i)\subseteq \overline{\omega_i}, \quad \sum_{i=1}^{N}\chi_{i} =1, }\\[4mm]
{\displaystyle \chi_{i}\in W^{1,\infty}(\omega_{i}),\;\;\Vert\nabla \chi_{i} \Vert_{L^{\infty}(\omega_i)} \leq \frac{C_{\chi}}{\mathrm{diam}\,(\omega_{i})}.}
\end{array}
\end{equation}
Note that the last inequality in \eqref{eq:2-0-1} implies that $\delta_{i}\geq c\,\mathrm{diam}\,(\omega_{i})$ for some constant $c>0$ independent of $i$, where $\delta_{i}$ denotes the `thickness' of the overlapping zone of $\omega_i$. This case is often referred to as \textit{generous overlap}, which we assume here and throughout the rest of this paper. Let $I_{h}:C(\Omega)\rightarrow U_{h}(\Omega)$ be the standard Lagrange interpolation, and \rs{define the following partition of unity operators} $P_{h,i}:  \big(U_{h,\mathcal{B}}(\omega^{\ast}_i), \Vert \cdot\Vert_{\mathcal{A},{\omega^{\ast}_{i}}}\big) \rightarrow \big(U_{h,0}(\omega^{\ast}_i), \Vert\cdot\Vert_{\mathcal{A},\omega^{\ast}_{i},k}\big)$ by
\begin{equation}\label{eq:5-3-0}
P_{h,i}v_{h} = I_{h}(\chi_{i}v_{h}).
\end{equation}
Here we equip the spaces $U_{h,\mathcal{B}}(\omega^{\ast}_i)$ and $U_{h,0}(\omega^{\ast}_i)$ with the norms $\Vert \cdot\Vert_{\mathcal{A},{\omega^{\ast}_{i}}}$ and $\Vert\cdot\Vert_{\mathcal{A},\omega^{\ast}_{i},k}$, respectively, and identify $U_{h,0}(\omega^{\ast}_i)$ as a subspace of $U_{h}(\Omega)$. On each $\omega_i^{\ast}$, we consider the eigenvalue problem of finding $(\phi_{h},\lambda_{h})\in U_{h,\mathcal{B}}(\omega_{i}^{\ast})\times \mathbb{R}$ such that
\begin{equation}\label{localOpeEVP}
P^{\ast}_{h,i}P_{h,i} \phi_{h} = \lambda_{h} \phi_{h},    
\end{equation}
where $P^{\ast}_{h,i}$ denotes the adjoint of $P_{h,i}$. The eigenproblem \eqref{localOpeEVP} has the following equivalent variational formulation:
\begin{equation}\label{localEVP}
\mathcal{A}_{\omega^{\ast}_{i},k}(I_{h}(\chi_{i} \phi_{h}), I_{h}(\chi_{i} v_{h})) = \lambda_{h}\,\mathcal{A}_{\omega_{i}^{\ast}}(\phi_{h}, v_{h})\quad \forall v_{h}\in U_{h,\mathcal{B}}(\omega_{i}^{\ast}).
\end{equation}
For each $j\in\mathbb{N}$, let $(\lambda_{h,j},\phi_{h,j})\in \mathbb{R}\times U_{h,\mathcal{B}}(\omega_{i}^{\ast})$ be the $j$-th eigenpair of problem \eqref{localEVP} (with eigenvalues sorted in decreasing order). The local approximation space $S_{n_i}(\omega_i)$ is defined by 
\begin{equation}\label{localappspace}
S_{n_i}(\omega_i) = {\rm span}\big\{\phi_{h,1}, \cdots, \phi_{h,n_i}\big\}.
\end{equation}
\begin{remark}\label{rem:n-widths}
Using the theory of the Kolmogrov $n$-width \cite{pinkus1985n,ma2023wavenumber}, it can be shown that $P_{h,i} (S_{n}(\omega_i))$ is the optimal $n$-dimensional space for the $n$-width
\begin{equation}\label{n-width}
d_{n}(P_{h,i}) := \inf_{Q(n)\subset U_{h,0}(\omega^{\ast}_{i})}\sup_{u_{h}\in U_{h,\mathcal{B}}(\omega^{\ast}_i)} \inf_{v_{h}\in Q(n)}\frac {\Vert P_{h,i}u_{h}-v_{h}\Vert_{\mathcal{A},\omega^{\ast}_{i},k}}{\Vert u_{h} \Vert_{\mathcal{A},\omega_{i}^{\ast}}},
\end{equation}
and that $d_{n}(P_{h,i})=\lambda^{1/2}_{h,n+1}$. 
\end{remark}

\begin{remark}
In the definition of eigenproblem \eqref{localEVP}, we can replace the bilinear form $\mathcal{A}_{\omega^{\ast}_{i}}(\cdot,\cdot)$ on the right-hand side with $\mathcal{A}_{\omega^{\ast}_{i},k}(\cdot,\cdot)$ (or equivalently, using the norm $\Vert\cdot\Vert_{\mathcal{A},\omega^{\ast}_{i},k}$ for $U_{h,\mathcal{B}}(\omega^{\ast}_i)$). This leads to the same eigenproblem as proposed in \cite{hu2024novel}. The modification does not change the theoretical results below, and makes little difference to numerical results.   
\end{remark}

With the local particular functions and local approximation spaces, we can define the global particular function and the global approximation space for MS-GFEM as follows. 
\begin{equation}\label{global_function_and_space}
 u_{h}^{p} :=\sum_{i=1}^{N}I_{h}(\chi_{i}\psi_{h,i})\;\;\;{\rm and}\;\; \; S_{n}(\Omega):= \Big\{ \sum_{i=1}^{N} I_{h}(\chi_{i}v_{h,i}):\; v_{h,i}\in S_{n_i}(\omega_{i}) \Big\}. 
\end{equation}
It is clear that $u_{h}^{p}\in U_{h}(\Omega)$ and $S_{h}(\Omega)\subset U_{h}(\Omega)$. Now we are ready to define the MS-GFEM approximation: Find $u^{G}_{h} = u_{h}^{p} + u_{h}^{s}$, with $u_{h}^{s}\in S_{n}(\Omega)$, such that
\begin{equation}\label{MSGFEM-approx}
\mathcal{B}(u_{h}^{G},\,v_{h}) = F(v_{h}) \quad \forall v_{h}\in S_{h}(\Omega).
\end{equation}
Using \eqref{fineFE_problem}, we see that \eqref{MSGFEM-approx} is equivalent to the following coarse-space correction:
\begin{equation}\label{MSGFEM-approx-equiv}
\mathcal{B}(u_{h}^{s},\,v_{h}) = \mathcal{B}(u^{{\mathdutchcal{e}}}_{h}-u_{h}^{p},\,v_{h}) \quad \forall v_{h}\in S_{h}(\Omega).
\end{equation}

As in Assumption~\ref{ass:global_stability}, we assume the following discrete inf-sup condition on the local sesquilinear form $B_{\omega_i^{\ast},{\rm imp}}(\cdot,\cdot)$; see also Remarks~\ref{rem:well-posedness} and \ref{rem:stab_constant}.
\begin{assumption}\label{ass:local_stability}
There exists an $h_{0}>0$ such that if $0<h< h_{0}$, the estimate
\begin{equation}
\inf_{v_{h}\in U_{h}(\omega_i^{\ast})} \sup_{0\neq u_{h}\in U_{h}(\omega_i^{\ast})} \frac{|B_{\omega_i^{\ast},{\rm imp}}(v_{h}, u_{h})|}{ \Vert v_{h}\Vert_{\mathcal{A},\omega^{\ast}_{i},k} \Vert u_{h}\Vert_{\mathcal{A},\omega^{\ast}_{i},k}} \gtrsim \frac{1 }{1 + kC_{{\mathtt{stab}},i}(k)},  
\end{equation}
holds with constant $C_{{\mathtt{stab}},i}(k)>0$ depending polynomially on $k$.
\end{assumption}

\begin{remark}\label{rem:local_infsup}
Following the proof of \cite[Proposition A.3]{chaumont2020multiscale}, it can be shown that if ${\rm diam}(\omega_i^{\ast})\sim k^{-1}$, then $C_{{\mathtt{stab}},i}(k)\sim k^{-1}$. Moreover, in the constant-coefficient case, if $\omega_i^{\ast}$ is star-shaped with respect to a ball of radius $H_i^{\ast}$, then $C_{{\mathtt{stab}},i}(k)\sim H_i^{\ast}$; see Remark~\ref{rem:well-posedness}.    
\end{remark}

We conclude this subsection by discussing the solution of the eigenproblem \eqref{localEVP} in practical computations. \medskip 

\noindent\textbf{Practical solution of eigenproblem \eqref{localEVP}}. The usual approach to solve \eqref{localEVP} is to first construct a basis for the Helmholtz-harmonic space $U_{h,\mathcal{B}}(\omega^{\ast}_i)$, and then solve the resulting matrix eigenvalue problem using standard eigensolvers. However, the first step is costly, especially when the mesh size $h$ is very small -- it requires solving $m_i$ local Helmholtz problems, with $m_i$ equal to the number of degrees of freedom on $\partial \omega_i^{\ast}$. In \cite{ma2023wavenumber}, it was suggested to solve an equivalent, augmented eigenproblem posed on standard finite element spaces. More precisely, we consider the following mixed formulation: find $\phi_h\in U_{h}(\omega_i^{\ast})$, $p_h\in U_{h,0}(\omega_i^{\ast})$, and $\lambda_h\in\mathbb{R}$ such that
\begin{equation}\label{mixed_EVP}
\begin{aligned}
\mathcal{A}_{\omega_{i}^{\ast}}(\phi_h,v_{h}) + \mathcal{B}_{\omega_{i}^{\ast}}(v_h,p_h)  =&\, \lambda^{-1}_{h}\mathcal{A}_{\omega^{\ast}_{i},k}\big(I_{h}(\chi_{i} \phi_{h}), I_{h}(\chi_{i} v_h)\big) \;\;\;\, \forall v_h\in U_{h}(\omega_{i}^{\ast}),\\
\mathcal{B}_{\omega_{i}^{\ast}}(\phi_{h},\xi_h) =&\;0\qquad \qquad\qquad\qquad\qquad\qquad\quad \;\forall \xi_{h}\in U_{h,0}(\omega_{i}^{\ast}).
\end{aligned}
\end{equation}
\rs{Although the problem size is about twice as large}, the eigenproblem \eqref{mixed_EVP} can be solved directly without constructing a basis for $U_{h,\mathcal{B}}(\omega^{\ast}_i)$. Note that by Assumption~\ref{ass:local_stability}, the sesquilinear form $\mathcal{B}_{\omega_{i}^{\ast}}(\cdot,\cdot)$ is inf-sup stable on $U_{h}(\omega_i^{\ast})\times U_{h,0}(\omega_i^{\ast})$. Hence, when writing \eqref{mixed_EVP} as a matrix eigenvalue problem, the matrix on the left-hand side is invertible; see \cite[Theorem 4.16]{ma2023wavenumber}. Furthermore, thanks to the exponential decay of the eigenvalues (see Theorem~\ref{thm:exponential-decay} below), \eqref{mixed_EVP} \RS{has large} (relative) spectral gap\rs{s between consecutive eigenvalues}. This advantageous property makes most iterative eigensolvers perform well in terms of convergence speed.

\subsection{Convergence theory of MS-GFEM}
In this subsection, we give convergence estimates for the MS-GFEM presented above. We start with local approximation error estimates. The following lemma is a direct consequence of Remark~\ref{rem:n-widths}; see also \cite[Theorem 2.10]{ma2023wavenumber}.
\begin{lemma}\label{lem:3-1}
Let the local particular function $\psi_{h,i}$ and the local approximation space $S_{n_i}(\omega_i)$ be defined by \eqref{localHelm_problem} and \eqref{localappspace}, respectively, and let $u^{{\mathdutchcal{e}}}_{h}$ be the solution of problem \eqref{fineFE_problem}. Then, 
\begin{equation}\label{local_error_bound}
\inf_{\varphi_{h}\in \psi_{h,i} + S_{n_{i}}(\omega_{i})}\big\Vert I_{h}\big(\chi_{i}(u^{\mathdutchcal{e}}_{h} - \varphi_{h})\big)\big\Vert_{\mathcal{A},\omega_{i},k}\leq \lambda^{1/2}_{h,n_i+1}\,\Vert u^{{\mathdutchcal{e}}}_{h}- \psi_{h,i}\Vert_{\mathcal{A},\omega_{i}^{\ast}},
\end{equation}
where $\lambda_{h,n_i+1}$ denotes the $(n_i+1)$-th eigenvalue of problem \eqref{localEVP}.
\end{lemma}

By Lemma~\ref{lem:3-1}, to prove an error bound for the local approximations, we need to estimate the two terms on the right-hand side of \eqref{local_error_bound}, i.e., $\Vert u^{{\mathdutchcal{e}}}_{h}- \psi_{h,i}\Vert_{\mathcal{A},\omega_{i}^{\ast}}$ and $\lambda^{1/2}_{h,n_i+1}$. We first consider the first term. In \cite{ma2023wavenumber}, it was bounded in terms of the source and boundary data $f$ and $g$, by means of stability estimates for the global and local Helmholtz problems. Here, in order to prove convergence for the iterative method defined below, we need to bound it in terms of a suitable norm of $u^{\mathdutchcal{e}}_{h}$. To this end, let us fix some notation. Let $\widetilde{\omega_i}^{\ast}$ denote the domain formed by enlarging $\omega_i^{\ast}$ with one layer of elements \rs{along the artificial boundary}, i.e.,
\begin{equation}\label{enlarged_domain}
\widetilde{\omega_i}^{\ast}: = \bigcup\big\{K\in \tau_{h}: \overline{K}\cap \overline{\omega_i^{\ast}} \neq \emptyset  \big\}.    
\end{equation}
Moreover, we set $H_{i}^{\ast}:= {\rm diam}(\omega^{\ast}_{i})$. The following lemma gives an estimate explicit in all parameters for $\Vert u^{{\mathdutchcal{e}}}_{h}- \psi_{h,i}\Vert_{\mathcal{A},\omega_{i}^{\ast}}$. It was proved in \cite{hu2024novel} for the classical Helmholtz equation, and the proof can be easily adapted to the heterogeneous-coefficient case by means of the discrete inf-sup condition in Assumption~\ref{ass:local_stability}.
\begin{lemma}\label{lem:error_local_Helmholtz}
Let Assumption~\ref{ass:local_stability} be satisfied with the stability constant $C_{{\mathtt{stab}},i}(k)$, and let $\widetilde{\omega_i}^{\ast}$ be defined by \eqref{enlarged_domain}. Then,
\begin{equation}\label{eq:error_local_Helmholtz}
\Vert u^{{\mathdutchcal{e}}}_{h}- \psi_{h,i}\Vert_{\mathcal{A},\omega_{i}^{\ast},k} \lesssim \big(1+kC_{{\mathtt{stab}},i}(k)\big) \big(1+ (kH_i^{\ast})^{-1}\big)^{\frac32} (kh)^{-\frac12}  \Vert u^{{\mathdutchcal{e}}}_{h}\Vert_{\mathcal{A},\widetilde{\omega_i}^{\ast},k}.
\end{equation} 
\end{lemma}
\begin{remark}\label{rem:kh-term}
Estimate \eqref{eq:error_local_Helmholtz} is pessimistic in terms of the dependence on the term $(kh)^{-\frac12}$. Indeed, numerical results in \cref{subsubsec:frequency_scaling} \alber{show} almost no dependence of the error on this term. 
\end{remark}
\begin{remark}\label{rem:DBC_local_stability}
If the local particular function $\psi_{h,i}$ is defined with a zero Dirichlet condition as in Remark~\ref{rem:DBC_localsolution}, then the following stability estimate holds when $H_{i}^{\ast} \lesssim k^{-1}$:
\begin{equation}\label{eq:improved_error_local_Helmholtz}
\Vert u^{{\mathdutchcal{e}}}_{h}- \psi_{h,i}\Vert_{\mathcal{A},\omega_{i}^{\ast},k} \lesssim  \Vert u^{{\mathdutchcal{e}}}_{h}\Vert_{\mathcal{A},\omega_{i}^{\ast},k}.  
\end{equation}
The proof is standard by means of Friedrichs's inequality.
\end{remark}

It remains to estimate $\lambda^{1/2}_{h,n_i+1}$, which lies in the core of the MS-GFEM theory. Combining \cite[Theorem 5.3]{ma2023wavenumber} and \cite[Theorem 3.8 and Theorem 6.19]{ma2023unified}, we have the following wavenumber-explicit exponential decay estimates for $\lambda^{1/2}_{h,n_i+1}$. For brevity, we drop the subdomain index $i$. 

\begin{theorem}\label{thm:exponential-decay}
Let $\delta^{\ast} = {\rm dist}(\omega,\,\partial \omega^{\ast}\setminus \partial \Omega)>0$, $H^{\ast} = {\rm diam}(\omega^{\ast})$, and let $d$ be the spatial dimension of $\Omega$. 
\begin{itemize}
    \item[(i)] There exist $n_{0}\sim (H^{\ast}/\delta^{\ast})^{d}$ and $b \sim (H^{\ast}/\delta^{\ast})^{-d/(d+1)}$, such that for any $n>n_{0}$, if $h\leq \min\{k^{-1},\delta^{\ast}/(bn^{1/(d+1)})\}$, then
\begin{equation}\label{eq:6-15}
\lambda^{1/2}_{h,n+1}\lesssim e^{\sigma}e^{-bn^{1/(d+1)}} \quad \text{with}\quad \sigma\sim k\delta^{\ast}.
\end{equation}
    \item[(ii)] There exist $n_{0}\sim (kH^{\ast})^{d^{2}/(d-1)} (H^{\ast}/\delta^{\ast})^{-1}$ and $b \sim (H^{\ast}/\delta^{\ast})^{-(d-1)/d}$, such that for any $n>n_{0}$, if $h\leq \min\{k^{-1},\delta^{\ast}/(bn^{1/d})\}$, then
\begin{equation}\label{eq:6-16}
\lambda^{1/2}_{h,n+1}\lesssim e^{-bn^{1/d}}.
\end{equation}
\end{itemize}
\vspace{-3ex}
\end{theorem}
\begin{remark}
The constants $n_0$, $b$, and $\sigma$ above do not depend on $k$, $h$, $H^{\ast}$, and $\delta^{\ast}$. For their precise definitions, see \cite[Theorem 5.3]{ma2023wavenumber} and \cite[Theorem 3.8 and Theorem 6.19]{ma2023unified}.      
\end{remark}
\begin{remark}
The estimates above show that the decay of $\lambda_{h,n+1}$ with respect to $n$ is independent of the wavenumber $k$ if $H^{\ast}\sim \delta^{\ast}\sim k^{-1}$ or if $n$ is sufficiently large. Moreover, the decay rate $O(e^{-bn^{1/d}})$ that we have proved is suboptimal compared with numerical results (see \cite{ma2023wavenumber}). Indeed, we conjecture that the actual decay rate is $O(e^{-bn^{1/(d-1)}})$.    
\end{remark}

Let $\xi$ (\RS{resp.} $\widetilde{\xi^{\ast}}$) be the maximum number of $\omega_i$'s (resp. $\widetilde{\omega_i^{\ast}}$'s) that overlap at any given point, i.e.,
\begin{equation}\label{coloring-constant}
\xi:= \max_{{\bm x} \in \Omega}{\Big( \text{card}\{i \ | \ {\bm x} \in \omega_i \} \Big)},\quad \widetilde{\xi^{\ast}}:= \max_{{\bm x} \in \Omega}{\Big( \text{card}\{i \ | \ {\bm x} \in \widetilde{\omega_i^{\ast}} \} \Big)}.
\end{equation}
Combining \rs{Lemma~\ref{lem:3-1}} and \rs{Lemma}~\ref{lem:error_local_Helmholtz}, and following the lines of the proof of \cite[Lemma 3.10]{ma2023wavenumber}, we get the following global approximation error estimate.
\begin{lemma}\label{lem:global_approximation_error}
Let the global particular function $u^{p}_{h}$ and the global approximation space $S_{h}(\Omega)$ be defined by \eqref{global_function_and_space}, and let $u^{{\mathdutchcal{e}}}_{h}$ be the solution of problem \eqref{fineFE_problem}. Then, 
\begin{equation}
\begin{array}{ll}
 {\displaystyle  \inf_{\varphi_{h}\in u^{p}_{h} + S_{h}(\Omega)} \frac{\big\Vert u^{\mathdutchcal{e}}_{h} - \varphi_{h}\big\Vert_{\mathcal{A},k}}{\big\Vert u^{\mathdutchcal{e}}_{h} \big\Vert_{\mathcal{A},k}} } \\[3ex]
{\displaystyle \quad \lesssim \sqrt{\xi \widetilde{\xi}^{\ast}} (kh)^{-\frac12} \max_{i=1,\cdots,N}\left(\big(1+kC_{{\mathtt{stab}},i}\big) \big(1+ (kH_i^{\ast})^{-1}\big)^{\frac32} \lambda^{1/2}_{h,n_i+1} \right).}
\end{array}
\end{equation}
\end{lemma}

We proceed to show the quasi-optimal convergence of the method. For this purpose, we need certain resolution conditions on the eigenvalue tolerance and on the size of the oversampling domains. For convenience, we set
\begin{equation}
 \lambda_{\mathtt{max}} = \max_{i=1,\cdots,N} \lambda_{h, n_i+1},\quad H^{\ast}_{\mathtt{max}}  =  \max_{i=1,\cdots,N} H^{\ast}_{i}.   
\end{equation}

\begin{theorem}\cite[Theorem 5.6]{ma2023wavenumber}\label{thm:quasi-optimality}
Let $u_{h}^{\mathdutchcal{e}}$ be the solution of problem \eqref{fineFE_problem} and $u_{h}^{G}$ be the MS-GFEM approximation defined by \eqref{MSGFEM-approx}. With $C_{{\mathtt{stab}}}(k)$ defined in Assumption~\ref{ass:global_stability}, if
\begin{equation}\label{resolution-conditions}
\lambda^{\frac12}_{\mathtt{max}} \lesssim \big(kC_{{\mathtt{stab}}}(k)\big)^{-1},\quad H^{\ast}_{\mathtt{max}} \lesssim k^{-1},
\end{equation}
then problem \eqref{MSGFEM-approx} is uniquely solvable with the estimate
\begin{equation}\label{eq:4-5-18}
\Vert u_{h}^{\mathdutchcal{e}}-u_{h}^{G}\Vert_{\mathcal{A},k} \leq 2C_{\mathcal{B}}\inf_{\varphi_{h}\in u_{h}^{p} + S_{h}(\Omega)}\Vert u_{h}^{\mathdutchcal{e}}-\varphi_{h} \Vert_{\mathcal{A},k},
\end{equation}
where $C_{\mathcal{B}}$ is defined in \eqref{boundedness-estimate}.
\end{theorem}
Combining \rs{Lemma}~\ref{lem:global_approximation_error} and \rs{Theorem}~\ref{thm:quasi-optimality} gives the following approximation result.
\begin{corollary}\label{cor:final_approx_result}
Let the assumptions in Theorem~\ref{thm:quasi-optimality} be satisfied. Then, 
\begin{equation}\label{eq:final_error_estimate}
\Vert u_{h}^{\mathdutchcal{e}}-u_{h}^{G}\Vert_{\mathcal{A},k} \leq \Lambda \,\Vert u^{\mathdutchcal{e}}_{h} \Vert_{\mathcal{A},k}  
\end{equation}
\rs{with}
\begin{equation}\label{convergence-rate}
 \Lambda \sim  \sqrt{\xi \widetilde{\xi}^{\ast}} C_{\mathcal{B}} (kh)^{-\frac12} \max_{i=1,\cdots,N}\left(\big(1+kC_{{\mathtt{stab}},i}(k)\big) \big(1+ (kH_i^{\ast})^{-1}\big)^{\frac32} \lambda^{1/2}_{h,n_i+1} \right).   
\end{equation}
\end{corollary}
\begin{remark}\label{rem:DBC_final_error_result}
The term $(kh)^{-\frac12}$ in \eqref{convergence-rate} arises from Lemma~\ref{lem:error_local_Helmholtz}. When using the local particular functions $\psi_{h,i}$ with a zero Dirichlet condition defined in Remark~\ref{rem:DBC_localsolution}, we have the improved estimate \eqref{eq:improved_error_local_Helmholtz} under the resolution conditions in \eqref{resolution-conditions}. Therefore, in this case, the error $\Lambda$ in \eqref{eq:final_error_estimate} satisfies 
\begin{equation}
\Lambda \lesssim \sqrt{\xi \widetilde{\xi}^{\ast}} C_{\mathcal{B}} \lambda^{\frac12}_{\mathtt{max}},   
\end{equation} 
without depending on the term $(kh)^{-\frac12}$. 
\end{remark}


We end this subsection by further discussing the error $\Lambda$ in Corollary~\ref{cor:final_approx_result}. Assume that 
\begin{equation}\label{new_resolution_conditions}
\lambda^{\frac12}_{h, n_i+1} \sim \big(kC_{{\mathtt{stab}}}(k)\big)^{-1}, \quad H^{\ast}_{i}\sim k^{-1}\qquad \text{for all}\;\; i=1,\cdots, N.  
\end{equation}
Then by Remark~\ref{rem:local_infsup}, we see that $kC_{{\mathtt{stab}},i}(k)\sim 1$ and thus
\begin{equation*}
\Lambda \sim (kh)^{-\frac12} \lambda^{\frac12}_{\mathtt{max}} \sim (kh)^{-\frac12} \big(kC_{{\mathtt{stab}}}(k)\big)^{-1}.    
\end{equation*}
It follows that $\Lambda\lesssim k^{-\alpha}$ for some $\alpha >0$, with a reasonable choice of the mesh-size $h$ (more precisely, $h\lesssim k^{-3+2\alpha}(C_{{\mathtt{stab}}}(k)\big)^{-2}$). Next we focus on the constant-coefficient case, and derive a more precise characterization for $\Lambda$. To do this, we make two reasonable assumptions.


\begin{assumption}\label{ass:mesh-size}
 The mesh size $h$ satisfies that $h\sim k^{-(1+\gamma)}$ for some $\gamma\in (0,1]$.  
\end{assumption}

\begin{assumption}\label{ass:star-shaped domain}
$\Omega$ is star-shaped with respect to a ball.   
\end{assumption}

Assumption~\ref{ass:mesh-size} is necessary for the relative error of the finite element approximation to remain bounded independent of $k$. In general, the choice of $\gamma$ depends on the order $p$ of the finite element approximation space. In particular, \cite{du2015preasymptotic} showed that $\gamma = \frac{1}{2p}$ is enough for a bounded (relative) error. Combining Remark~\ref{rem:well-posedness} and Assumption~\ref{ass:star-shaped domain}, we see that $C_{\mathtt{stab}}(k)\sim 1$. Under these two assumptions and the conditions in \eqref{new_resolution_conditions}, we deduce that
\begin{equation}\label{asymptotic_decay}
    \Lambda \sim k^{\frac{\gamma}{2}} \lambda^{\frac12}_{\mathtt{max}} \sim   (\lambda^{\frac12}_{\mathtt{max}})^{1-\frac{\gamma}{2}} \sim k^{-1+\frac{\gamma}{2}}.
\end{equation}
Namely, the error of MS-GFEM decays with increasing $k$ as $k^{-1+\frac{\gamma}{2}}$. Moreover, as noted in Remark~\ref{rem:kh-term}, the result above might be improved in terms of the dependence on the term $k^{\frac{\gamma}{2}}$.

\section{MS-GFEM based two-level RAS methods}\label{sec:preconditioner}
\subsection{Iterative MS-GFEM}
Before defining the iterative \alber{MS-GFEM}, let us briefly describe the idea behind it. Note that in the MS-GFEM theory, we only need the local particular functions $\psi_{h,i}$ to satisfy the underlying PDE locally, and they can have arbitrary boundary conditions (as long as the local problems are solvable). In Section~\ref{sec:MS-GFEM}, we have \rs{made two explicit choices and} used homogeneous impedance/Dirichlet boundary conditions for $\psi_{h,i}$. Indeed, if more ‘correct' boundary conditions are available and used for $\psi_{h,i}$, then $\psi_{h,i}$ will be a better approximation of $u_{h}^{e}|_{\omega_i^{\ast}}$, and thus a more accurate (whole) local approximation can be obtained (see Lemma~\ref{lem:3-1}). Given a MS-GFEM approximation $u_{h}^{j}$, a natural idea to obtain a more accurate approximation is to update $\psi_{h,i}$ using the local boundary conditions of $u_{h}^{j}$. In the rest of this section, we focus on the MS-GFEM with impedance conditions for the $\psi_{h,i}$'s because of its superior performance, but the proposed preconditioner and its convergence estimates also hold for the case of local Dirichlet conditions.    


To formulate the method \RS{more} precisely, let us introduce some notation. Recall the local FE spaces $U_{h}(\omega_i^*)$, the coarse approximation space $S_{n}(\Omega)$, and the nodewise extension operators $E_{i}$. Throughout this section, we fix the choices \rs{of }$\omega_i$, $\omega_i^*$, $\chi_i$, and $S_{n}(\Omega)$. For each $i=1,\cdots,N$, we define \rs{a} local projection operator $\pi_i: U_{h}(\Omega) \rightarrow U_{h}(\omega_i^{\ast})$ such that for each $v_{h}\in U_{h}(\Omega)$, $\pi_i(v_{h})\in U_{h}(\omega_i^*)$ satisfies
\begin{equation}\label{local_projection}
 \mathcal{B}_{\omega_i^{\ast},{\rm imp}}(\pi_i(v_h),w_h) = \mathcal{B}(v_h,E_i w_h) \quad \forall w_h \in U_{h}(\omega_i^*).
\end{equation}
Similarly, we define the coarse projection operator $\pi_S:  U_{h}(\Omega)\rightarrow S_{n}(\Omega)$ \rs{such that} for each $v_{h}\in U_{h}(\Omega)$, $\pi_S(v_{h})\in S_{n}(\Omega)$ satisfies
\begin{equation}\label{coarse_projection}
\mathcal{B}(\pi_S(v_h),w_h) = \mathcal{B}(v_h,w_h) \quad \forall w_h \in S_{n}(\Omega).       
\end{equation}
Here we assume that Assumption~\ref{ass:local_stability} and the conditions in Theorem~\ref{thm:quasi-optimality} are satisfied such that the operators $\pi_i$ and $\pi_S$ are well-defined. For each $i=1,\cdots,N$, we also define the operator ${\chi}_{h,i}: U_{h}(\omega_i^*)\rightarrow U_{h}(\Omega)$ by
$${\chi}_{h,i}(v_{h}) = I_{h}(\chi_{i} v_{h}),$$
where $\chi_i$ is the partition of unity function, and $I_{h}$ is the standard nodal interpolation. \rs{Using this notation}, we can \rs{now} define the iterative method following the idea described above. Given a MS-GFEM approximation $u_{h}^{j}$, we first define the (new) local particular functions 
\begin{equation}\label{updated_localup}
 \psi^{j+1}_{h,i}: =  u_{h}^{j}|_{\omega_i^{\ast}} + \pi_{i}\big(u_{h}^{e} - u_{h}^{j}\big)\in \RS{U_{h}(\omega_i^{\ast})},  \quad j=1,\cdots,N,
\end{equation}
where $u_{h}^{e}$ is the solution of problem \eqref{fineFE_problem}. The new MS-GFEM approximation $u_{h}^{j+1}$ is then defined similarly as in Section~\ref{sec:MS-GFEM} by first gluing the local particular functions together and then adding a coarse-space correction (see \eqref{MSGFEM-approx-equiv}):
\begin{equation}\label{iterative-MSGFEM-v1}
u_{h}^{j+1}= \sum\limits_{i=1}^N {\chi}_{h,i} (\psi^{j+1}_{h,i}) + \pi_S \Big(u_{h}^{e}  - \sum\limits_{i=1}^N {\chi}_{h,i}(\psi^{j+1}_{h,i}) \Big) .   
\end{equation}
\begin{remark}
In \cite{gong2023convergence}, it was shown that \eqref{updated_localup} is the finite element analogue of the following definition of $\psi^{j+1}_{i}$ at the PDE level: given $u^{j}$, find $\psi^{j+1}_{i}$ such that
\begin{equation}
\left\{
\begin{array}{lll}
{\displaystyle -{\rm div}(A\nabla \psi^{j+1}_{i}) - k^{2}V^{2}\psi^{j+1}_{i}= f\,\qquad\qquad\qquad\qquad{\rm in}\;\, \omega_i^{\ast}, }\\[3mm]
{\displaystyle A\nabla \psi^{j+1}_{i}\cdot {\bm n} - {\rm i}kV\psi^{j+1}_{i}=A\nabla u^{j}\cdot {\bm n} - {\rm i}kV u^{j} \qquad \quad \,{\rm on}\;\,\partial \omega_i^{\ast} \setminus \partial \Omega,}\\[3mm]
{\displaystyle A\nabla \psi^{j+1}_{i}\cdot {\bm n} - {\rm i}k\beta \psi^{j+1}_{i}=g  \qquad \qquad\qquad\qquad\quad\;\; \quad {\rm on}\;\,\partial \omega_i^{\ast}\cap \partial \Omega.}
\end{array}
\right.
\end{equation}
\RS{This shows clearly} that \RS{in fact we are using} $u^{j}$ as impedance boundary data for $\psi^{j+1}_{i}$. 
\end{remark}

\begin{remark}
Without the coarse-space correction, the iterative method \eqref{iterative-MSGFEM-v1} reduces to the Optimized Restricted Additive Schwarz (ORAS) method. Its convergence was proved for the classical Helmholtz equation in strip-shaped domains in \cite{gong2023convergence}.        
\end{remark}

In order to prove the convergence of the iterative method, we give an equivalent definition of \eqref{iterative-MSGFEM-v1} following \cite{strehlow2024fast}. We first define the following MS-GFEM map.
\begin{definition}
With the operators $\pi_i$ and $\pi_S$ defined in \eqref{local_projection} and \eqref{coarse_projection}, we define the map $G:U_{h}(\Omega) \rightarrow U_{h}(\Omega)$ by
    \begin{equation}\label{MSGFEM-map}
               G(v_h) = \sum\limits_{i=1}^N {\chi}_{h,i} \pi_i(v_h) + \pi_S \Big( v_h - \sum\limits_{i=1}^N {\chi}_{h,i} \pi_i(v_h) \Big). 
    \end{equation}
We call $G$ the MS-GFEM map.  
\end{definition}
It follows from the definition of the operators $\pi_i$ and $\pi_S$ that, for any $v_{h}\in U_{h}(\Omega)$, $G(v_{h})$ is exactly the MS-GFEM \rs{approximation for} the FE problem \eqref{fineFE_problem} with the right-hand side $F(\cdot):=\mathcal{B}(v_h, \cdot)$ (and thus with \rs{exact} solution $v_h$). Using this fact, we have the following estimate due to Corollary~\ref{cor:final_approx_result}.
\begin{lemma} \label{lem:convergence-mapG}
 Let $G$ be defined above, and let the assumptions in Theorem~\ref{thm:quasi-optimality} be satisfied. Then, for any $v_h \in U_{h}(\Omega)$,
 \begin{equation*}
\Vert v_h - G(v_h) \Vert_{\mathcal{A},k} \leq \Lambda \Vert v_{\alber{h}} \Vert_{\mathcal{A},k},    
 \end{equation*}
 where $\Lambda$ is given by \eqref{convergence-rate}.
  \end{lemma}
Now we can give an equivalent definition of \eqref{iterative-MSGFEM-v1} using the map $G$.
\begin{definition}[Iterative MS-GFEM] \label{iteration}
Let $u_{h}^{e}$ be the solution of problem \eqref{fineFE_problem}. Given the $j$-th iterate $u_{h}^{j} \in U_{h}(\Omega)$, set
    \begin{equation} \label{iterative-MSGFEM-v2}
     u_{h}^{j+1}:= u_{h}^{j} + G(u_{h}^{e} -u_{h}^{j}).   
    \end{equation}
\end{definition}
The equivalence of the two definitions is an easy consequence of the fact \rs{that}
\begin{equation*}
\sum\limits_{i=1}^N {\chi}_{h,i}(v_h|_{\omega_i^{\ast}}) = v_h \quad \forall v_h\in U_{h}(\Omega).    
\end{equation*}
We have the following convergence estimate for the iterative method as a direct result of Lemma~\ref{lem:convergence-mapG}.
\begin{proposition} \label{prop:conv-MSGFEM-iteration}
Let the iterates $\{u_{h}^{j}\}_{j \in \mathbb{N}}$ be defined by \eqref{iterative-MSGFEM-v2}, and let the assumptions in Theorem~\ref{thm:quasi-optimality} be satisfied. Then, 
    \begin{equation*}
   \Vert u_{h}^{e}- u_{h}^{j+1} \Vert_{\mathcal{A},k} \leq \Lambda \, \Vert u_{h}^{e}- u_{h}^{j} \Vert_{\mathcal{A},k}.  
    \end{equation*}
\end{proposition}


\subsection{Two-level RAS preconditioner}
The iterative MS-GFEM is essentially a preconditioned Richardson iterative method. To see this and to derive the preconditioner, we first write the method in matrix form. For each $i = 1,\cdots,N$, let $\mathbf{E}_{i}$ and ${\bm \chi}_i$ be the matrix representations of the operators $E_i$ and ${\chi}_{h,i}$ with respect to the basis $\{\phi_i\}_{i=1}^{m}$ of $U_{h}(\Omega)$ and a basis of $U_{h}(\omega_i^{\ast})$, respectively. Similarly, we let ${\bf E}_{0}$ be the matrix representation of the natural embedding $E_{0}: S_{n}(\Omega)\rightarrow U_{h}(\Omega)$ with respect to a chosen basis of $S_{n}(\Omega)$. Moreover, for each $i=1,\cdots,N$, we let $\mathbf{B}_i$ be the matrix such that for all $u_h, v_h\in U_{h}(\omega_i^{\ast})$,
\begin{equation}
\mathcal{B}_{\omega_i^{\ast},{\rm imp}} (v_h, u_h) = {\bf u}^{*} \mathbf{B}_i {\bf v},    
\end{equation}
where ${\bf u}$ and ${\bf v}$ are the vector representations of $u_h$ and $v_h$ with respect to the chosen basis of $U_{h}(\omega_i^{\ast})$, respectively, and ${\bf u}^{*}$ denotes the conjugate transpose of ${\bf u}$. \rs{Thus}, we can represent the operators $\pi_i$ and $\pi_S$ in matrix form:
$$
    \boldsymbol{\pi}_i = {\bf B}_i^{-1} {\bf E}^{\top}_i {\bf B} \quad \text{and} \quad \boldsymbol{\pi}_S = {\bf B}_0^{-1}{\bf E}^{\top}_0 {\bf B},
$$
where ${\bf E}^{\top}_i$ denotes the transpose of ${\bf E}_i$, ${\bf B}$ is the matrix of the linear system \eqref{linear_system}, and ${\bf B}_0:= {\bf E}^{\top}_0 {\bf B} {\bf E}_0$. As a result, the MS-GFEM map $G$ defined by \eqref{MSGFEM-map} \rs{in matrix form is given by}
\begin{align*}
  {\bf G} 
    &= \Bigg(  \sum_{i=1}^{N} {\boldsymbol{\chi}}_i {\bf B}_i^{-1} {\bf E}^{\top}_i {\bf B} \Bigg)
    + 
    \big({\bf E}_0 {\bf B}_0^{-1} {\bf E}^{\top}_0 {\bf B} \big) 
    \left( {\bf I} - \sum_{i=1}^{N} {\boldsymbol{\chi}}_i {\bf B}_i^{-1} {\bf E}^{\top}_i {\bf B} \right).
\end{align*}
Define
\begin{equation}\label{preconditioner}
    {\bf M}^{-1} := \Bigg(  \sum_{i=1}^{N} {\boldsymbol{\chi}}_i {\bf B}_i^{-1} {\bf E}^{\top}_i  \Bigg)
    + 
    \big({\bf E}_0 {\bf B}_0^{-1} {\bf E}^{\top}_0  \big) 
    \left( {\bf I} - {\bf B}\sum_{i=1}^{N} {\boldsymbol{\chi}}_i {\bf B}_i^{-1} {\bf E}^{\top}_i \right).
\end{equation}
It follows that ${\bf G} = {\bf M}^{-1} {\bf B}$. Recalling the linear system \eqref{linear_system}, \eqref{iterative-MSGFEM-v2} can then be written in matrix form as
\begin{equation}\label{iterative-MSGFEM-v3}
 {\bf v}_{j+1}: = {\bf v}_{j} + {\bf M}^{-1} {\bf B}({\bf u} - {\bf v}_{j}) = {\bf v}_{j} + {\bf M}^{-1} ({\bf F} - {\bf B}{\bf v}_{j}).
\end{equation}
Hence, \eqref{iterative-MSGFEM-v3} is a preconditioned Richardson iterative method for \eqref{linear_system}, with the two-level, hybrid RAS preconditioner ${\bf M}^{-1}$. Here ‘hybrid' means that the coarse space is incorporated multiplicatively. The preconditioner corresponds to the two-level Hybrid-II Schwarz preconditioner \cite[p. 48]{smith2004domain}, and is closely related to the ‘Adapted Deflation Variant 2' preconditioner proposed in \cite{tang2009comparison}. Moreover, we note that except for oversampling, our preconditioner corresponds to the one proposed in \cite{hu2024novel} for the classical Helmholtz equation. However, here we derive it naturally from MS-GFEM, as opposed to \cite{hu2024novel} where the \rs{coarse space} was \rs{added explicitly}.

\subsection{Convergence of preconditioned GMRES}
With ${\bf M}^{-1}$ defined in \eqref{preconditioner}, we consider using GMRES to solve the preconditioned linear system:
\begin{equation}\label{preconditioned_linear_system}
    {\bf M}^{-1} {\bf Bu} = {\bf M}^{-1} {\bf F}.
\end{equation}
Before proceeding, let us briefly describe the GMRES algorithm. Let $b(\cdot,\cdot)$ be an inner product on $\mathbb{C}^{m}$ and $\Vert\cdot\Vert_{b}$ the corresponding norm. Given an initial approximation ${\bf u}_{0}\in \mathbb{C}^{m}$ with the initial residual ${\bf r}_{0}= {\bf M}^{-1} ({\bf F}- {\bf Bu}_{0})$, in the $j$-th iteration, the GMRES algorithm applied to \eqref{preconditioned_linear_system} computes a correction vector ${\bf z}_j$ from the Krylov subspace
$$\mathcal{K}_j:= \text{span} \big\{ {\bf r}_{0},\,({\bf M}^{-1}{\bf B}){\bf r}_{0},...,({\bf M}^{-1}{\bf B})^{j-1}{\bf r}_{0} \big\},$$
such that ${\bf z}_j$ minimizes the $\Vert\cdot\Vert_{b}$-norm of the residual, i.e., ${\bf z}_j$ solves 
\begin{equation}\label{GMRES-minimization}
\min_{{\bf z} \in \mathcal{K}_j} \big\lVert {\bf M}^{-1} {\bf F} - {\bf M}^{-1}{\bf B}({\bf u}_{0} + {\bf z}) \big\rVert_b    
\end{equation}
The $j$-th iterate is then defined by ${\bf u}_j:= {\bf u}_{0} + {\bf z}_j$. 

The minimal residual property \eqref{GMRES-minimization} of GMRES has an important implication. Let $\{{\bf v}_j \}_{j\in\mathbb{N}}$ and $\{{\bf u}_j \}_{j\in\mathbb{N}}$ be the sequences of iterates generated by the MS-GFEM iteration \eqref{iterative-MSGFEM-v3} and the GMRES algorithm with \rs{preconditioner $M^{-1}$ and} the same initial guess ${\bf u}_{0} = {\bf v}_{0}$, respectively. Using a simple induction argument, we can easily show that 
$${\bf v}_{j}\in {\bf v}_{0}+\mathcal{K}_{j} = {\bf u}_{0}+\mathcal{K}_{j}.$$ 
Then it follows from the definition of the GMRES algorithm that 
\begin{equation}\label{residual_comparison}
\big\lVert {\bf M}^{-1} {\bf F} - {\bf M}^{-1}{\bf B}{\bf u}_{j} \big\rVert_b \leq \big\lVert {\bf M}^{-1} {\bf F} - {\bf M}^{-1}{\bf B}{\bf v}_{j} \big\rVert_b.    
\end{equation}
In other words, with residual measured in the $\Vert\cdot\Vert_{b}$-norm, the GMRES algorithm converges at least as fast as the MS-GFEM iteration. With an abuse of notation, we denote by $\Vert\cdot\Vert_{\mathcal{A},{k}}$ the norm on $\mathbb{C}^{m}$ induced by the $\Vert\cdot\Vert_{\mathcal{A},{k}}$-norm on $U_{h}(\Omega)$ (defined in \eqref{localforms}), i.e., for any ${\bf v}\in \mathbb{C}^{m}$,
\begin{equation}
 \Vert {\bf v}\Vert_{\mathcal{A},{k}}:= \Vert v_h\Vert_{\mathcal{A},{k}},
\end{equation}
where $v_{h}\in U_{h}(\Omega)$ is the \rs{finite element} function \rs{corresponding to the coefficient vector} ${\bf v}$ with respect to the basis $\{\phi_i\}_{i=1}^{m}$. To derive the rate of convergence for the GMRES algorithm, we make the following assumption concerning the equivalence of the two norms $\Vert\cdot\Vert_{b}$ and $\Vert\cdot\Vert_{\mathcal{A},{k}}$.

\begin{assumption}\label{ass:norm-equiv}
     There exist constants $b_1, b_2 > 0$ such that for all ${\bf v} \in \mathbb{C}^m$,
    $$
        b_1 \lVert {\bf v} \rVert _b \leq \lVert {\bf v} \rVert _{\mathcal{A},k} \leq b_2 \lVert {\bf v} \rVert _b.
    $$        
\end{assumption}

The following theorem gives an upper bound on the residual of GMRES applied to \eqref{preconditioned_linear_system}. The proof is similar to that of \cite[Theorem 3.9]{strehlow2024fast} for positive definite problems, based on the convergence estimate for the iterative MS-GFEM (Proposition~\ref{prop:conv-MSGFEM-iteration}) and \eqref{residual_comparison}.
\begin{theorem}\label{thm:gmresconvergence}
Let the assumptions in Corollary~\ref{cor:final_approx_result} be satisfied, and let $\Lambda<1$ be defined by \eqref{convergence-rate}. Then,
\begin{equation*}
  \big\lVert {\bf M}^{-1}{\bf B} \big\rVert_{\mathcal{A},k} \big\lVert \big({\bf M}^{-1}{\bf B}\big)^{-1} \big\rVert_{\mathcal{A},k} \leq \frac{1+\Lambda}{1-\Lambda}. 
\end{equation*}
Moreover, after $j$ steps, the residual ${\bf r}_{j}:= {\bf M}^{-1}{\bf F} - {\bf M}^{-1}{\bf B}{\bf u}_{j}$ of GMRES applied to \eqref{preconditioned_linear_system} satisfies 
    $$
        \lVert {\bf r}_{j} \rVert_b \leq \Lambda^j \left( \frac{1+\Lambda}{1-\Lambda}\right)\frac{b_2}{b_1} \lVert {\bf r}_{0} \rVert_b,\quad j=1,2\cdots,
    $$
where ${\bf r}_{0}$ is the initial residual, and $b_1, b_2>0$ are given in Assumption~\ref{ass:norm-equiv}. 
\end{theorem}
\begin{remark}
For ‘standard GMRES' with the Euclidean norm, \cite[Corollary 5.8]{gong2021domain} showed that $b_2/b_1 \sim (kh)^{-1}$ when the mesh family $\{\tau_h\}$ is quasiuniform. In this case, it needs at most extra $|\log_2(kh)|/|\log_2(\Lambda)|$ iterations to achieve the same relative residual as GMRES applied in the $\Vert\cdot\Vert_{\mathcal{A},k}$-norm. In the constant-coefficient case, as discussed at the end of Section~\ref{sec:MS-GFEM}, if $h\sim k^{-1-\frac{1}{2p}}$ ($p\geq 1$), then $\Lambda\sim k^{-1+\frac{1}{4p}}$. In this case, only one extra iteration is needed.
\end{remark}

Theorem~\ref{thm:gmresconvergence} indicates that GMRES applied to \eqref{preconditioned_linear_system} converges at least at a rate of $\Lambda$, with $\Lambda$ being the error of the underlying MS-GFEM. Notably, this result is not proved by means of the commonly used ‘Elman theory' \cite{eisenstat1983variational}, which requires bounds on the norm and fields of value of the preconditioner. More precisely, to use this theory in our setting, one needs to prove that there exist constants $c_1,c_2>0$ such that with ${\bf P}={\bf M}^{-1}{\bf B}$,
\begin{equation}\label{field of values}
\Vert {\bf P} {\bf v}\Vert_{b} \leq c_1 \Vert {\bf v}\Vert_{b}, \quad \text{and}\quad  \big|({\bf v}, {\bf P} {\bf v})_{b} \big|  \geq c_2 \Vert {\bf v}\Vert^{2}_{b}\qquad \text{for all}\;\; {\bf v}\in \mathbb{C}^{m}.  
\end{equation}
With the bounds in \eqref{field of values}, \cite{eisenstat1983variational} showed that after $j$ steps, the norm of the residual of GMRES is bounded by
\begin{equation}
\lVert {\bf r}_{j} \rVert_b \leq C\Big(1-\frac{c_2^{2}}{c_1^{2}}\Big)^{\frac{j}{2}} \lVert {\bf r}_{0} \rVert_b. 
\end{equation}
Unfortunately, the bound for the convergence rate given by this theory is far from sharp. In \cite{hu2024novel}, the authors established a $k$-independent convergence rate for GMRES by proving \eqref{field of values} with $k$-independent constants $c_1$ and $c_2$. In contrast, under the same conditions as \cite{hu2024novel}, we obtain a much sharper estimate for the convergence rate of GMRES in terms of the error of the underlying MS-GFEM that has been well analysed. In particular, the discussion at the end of Section~\ref{sec:MS-GFEM} shows that in the constant-coefficient case, the convergence rate $\Lambda$ decays with the eigenvalue tolerance $\lambda^{\frac12}_{\mathtt{max}}$ as $(\lambda^{\frac12}_{\mathtt{max}})^{1-\frac{\gamma}{2}}$, or equivalently, it decays with wavenumber $k$ as $k^{-1+\frac{\gamma}{2}}$.


\section{Numerical Experiments}\label{sec:numerical results}
In this section, we provide numerical results for two and three dimensional Helmholtz problems. We implement the two-dimensional experiments using FEniCS \cite{alnaes2015fenics}, and the three-dimensional experiments using FreeFEM \cite{hecht2012new}, in particular on top of its \texttt{ffddm} framework \cite{FFD:Tournier:2019}. The local eigensolves are performed using ARPACK, and the local Helmholtz solves and the coarse solve are performed using direct solvers like MUMPS. For all the tests, the stopping criterion of GMRES is a relative residual reduction of $10^{-6}$.

\subsection{Constant-coefficient Example}
In this example, we consider equation \eqref{eq:1-1} on the unit square $\Omega=(0,1)^2$ with $A({\bm x}) =I$, $V({\bm x}) = 1$, and $\beta({\bm x}) =1$. The underlying fine FE mesh with mesh-size $h$ on $\Omega$ is based on a uniform Cartesian grid. For the discretization, we use either linear or quadratic elements and choose the fine mesh-size $h\sim k^{-1-\frac{1}{2p}}$ ($p$ denotes the polynomial degree). To implement the MS-GFEM, we first partition the domain into $M=m^{2}$ uniform square, non-overlapping subdomains $\{\mathring{\omega}_i\} $ resolved by the mesh, and then extend each $\mathring{\omega}_i$ by adding one or several layers of adjacent fine mesh elements to create an overlapping decomposition $\{ \omega_{i}\}_{i=1}^{M}$ of $\Omega$. Each subdomain $\omega_{i}$ is further extended by adding several layers of fine mesh elements to create an oversampling domain $\omega_{i}^{\ast}$ on which the local problems are solved. We use $\delta: = {\rm dist}(\mathring{\omega}_i,\partial \omega_i\setminus \partial \Omega)$ and $\delta^{\ast}: = {\rm dist}(\omega_i,\partial \omega^{\ast}_i\setminus \partial \Omega)$ to denote the overlap and oversampling sizes, respectively. For illustration purposes, we also use the quantity $H/H^{\ast}$ as a measure of oversampling, where $H$ and $H^{\ast}$ denote the side lengths of the (interior) subdomains and the oversampling domains, respectively. Clearly, $H/H^{\ast} = 1$ means no oversampling is used. The number of eigenvectors used per subdomain is denoted by $n_{loc}$. 

The partition of unity functions $\{\chi_i\}_{i=1}^{M}$ are constructed as follows. As shown in \cref{fig:partition of unity}, the function $\chi_i$ over the subdomain $\omega_i$ is 1 on the non-overlapping zone $\mathring{\omega}_i$ (on the center in white), linear on the four strips of the overlap (in green), and bilinear on the four corners (in blue). 

\begin{figure}
\centering
\begin{tikzpicture}
\filldraw[thick, dotted, fill=blue!20] (0.8,0.8) rectangle (1.2,1.2);
\filldraw[thick, dotted,fill=green!20] (0.8,1.2) rectangle (1.2,1.8);
\filldraw[thick, dotted,fill=blue!20] (0.8,1.8) rectangle (1.2,2.2);
\filldraw[thick, dotted,fill=green!20] (1.2,1.8) rectangle (1.8,2.2);
\filldraw[thick, dotted,fill=blue!20] (1.8,1.8) rectangle (2.2,2.2);
\filldraw[thick, dotted,fill=green!20] (1.8,1.2) rectangle (2.2,1.8);
\filldraw[thick, dotted,fill=blue!20] (1.8,0.8) rectangle (2.2,1.2);
\filldraw[thick, dotted,fill=green!20] (1.2,0.8) rectangle (1.8,1.2);

\draw[->] (2.1,1.6) -- (3.7,1.0);
\node at (4.1, 0.9) {$\omega_i$};
\draw[thick] (0,0) grid (3,3);
\end{tikzpicture}
    \caption{Illustration of a subdomain $\omega_i$ with the overlapping zone shown in color. }
    \label{fig:partition of unity}
\end{figure}


\subsubsection{Frequency Scaling Test}\label{subsubsec:frequency_scaling}
\begin{figure}
    \centering
    \includegraphics[scale=0.4]{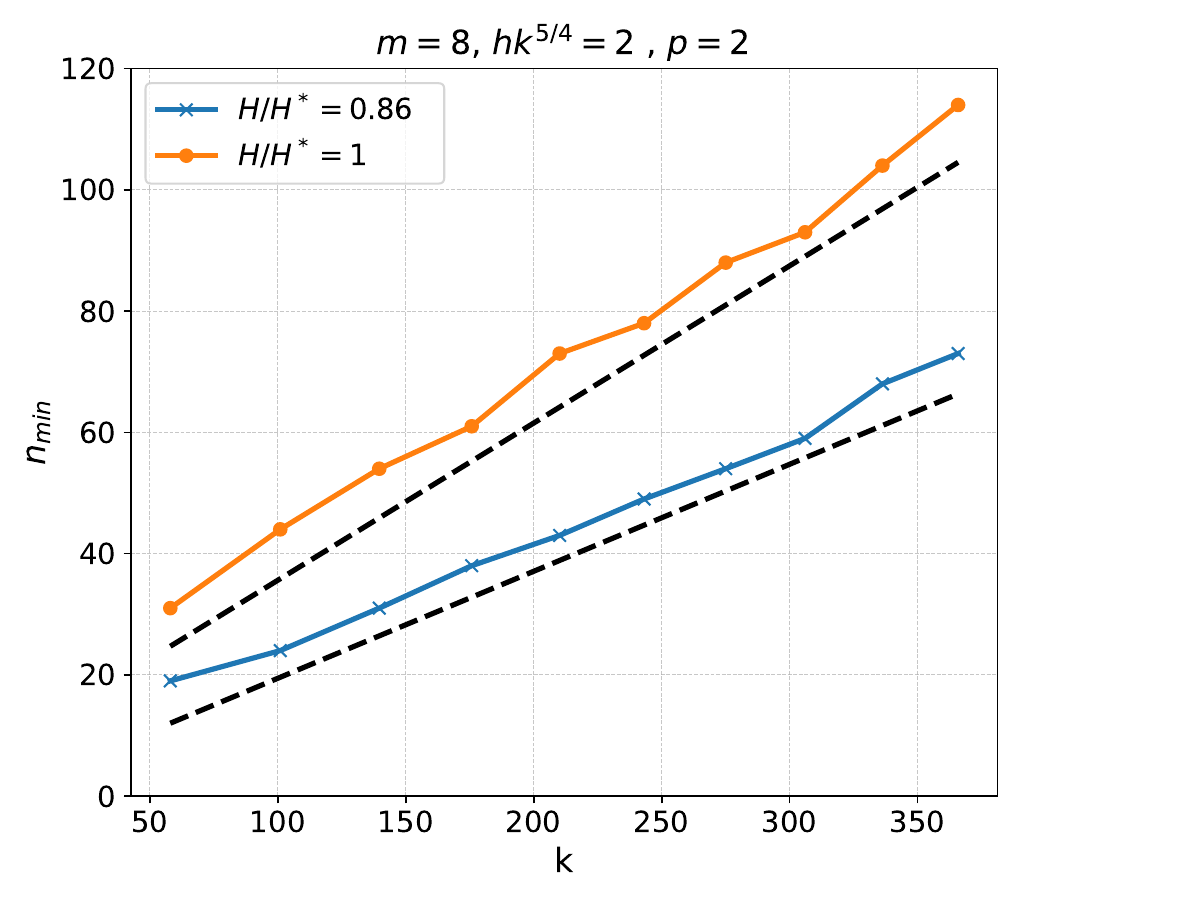}
    \includegraphics[scale=0.4]{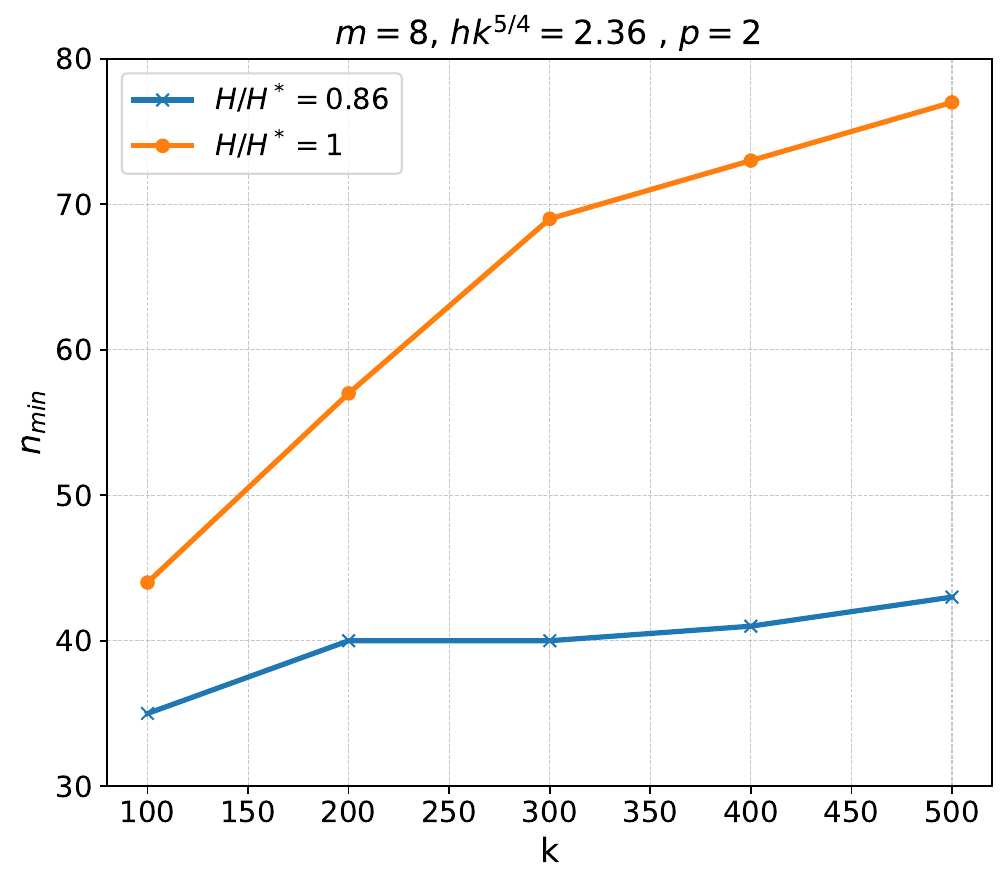}
    \caption{Frequency Scaling test:  the numbers of eigenfunctions $n_{min}$ needed to achieve the accuracy criterion \eqref{eq:accuracy_criterion}for case (i) (left, with $H$ fixed) and case (ii) (right, with $kH$ fixed).}
    \label{fig:scaling}
\end{figure} 

As a first step, we evaluate the approximation capability of the local approximation spaces for high-frequency waves. For this purpose, we introduce the following accuracy criterion for the MS-GFEM approximation: given a specific setup of the MS-GFEM, we define $n_{min}$ as the smallest number of eigenfunctions (we consider the same number of eigenfunctions on each subdomain) needed such that the MS-GFEM approximation $u_{h}^{G}$ satisfies
\begin{equation}
\label{eq:accuracy_criterion}
    |u_h^{\mathdutchcal{e}} - u_{h}^{G}|_{H^1(\Omega)} \leq |u^{\mathdutchcal{e}} - u_{h}^{\mathdutchcal{e}}|_{H^1(\Omega)},
\end{equation}
where $u^{\mathdutchcal{e}}$ and $u^{\mathdutchcal{e}}_{h}$ denote the exact solution and the standard FE approximation of the problem, respectively. We want to investigate how $n_{min}$ grows with increasing wavenumber $k$.

For the test, we choose the problem that admits the plane-wave solution $u^{\mathdutchcal{e}} = \exp\big({\rm i}\bm{k} \cdot {\bm x}\big)$ with $\bm{k} = k(0.6,\,0.8)$, and use P2 finite elements for the discretization. Two cases are considered: (i) the subdomains are fixed ($H\approx 12\delta \approx 1/8$) for increasing $k$ and (ii) the subdomains shrink in size with $H \approx 13\delta \approx 30k^{-1}$. The results are shown in \cref{fig:scaling} for the two cases with different oversampling ratios. We can see that for case (i), the number $n_{min}$ increases linearly with $k$, and it is significantly smaller if oversampling is used. These results show that the MS-GFEM approximation with $O(k)$ DOFs is comparable to the fine-scale FE approximation with $O(k^{\frac{5}{2}})$ DOFs in terms of accuracy -- albeit with a more densely populated linear system. Compared with case (i), we observe that $n_{min}$ in case (ii) grows much more slowly. Moreover, the effect of oversampling on $n_{min}$ is more prominent in the second case -- the growth of $n_{min}$ with $k$ is very mild with oversampling. Indeed, our theory predicts that $n_{min}$ grows at most like $O(\log^{2}(k))$ in this case.  



\begin{figure}
    \centering
    \includegraphics[scale=0.43]{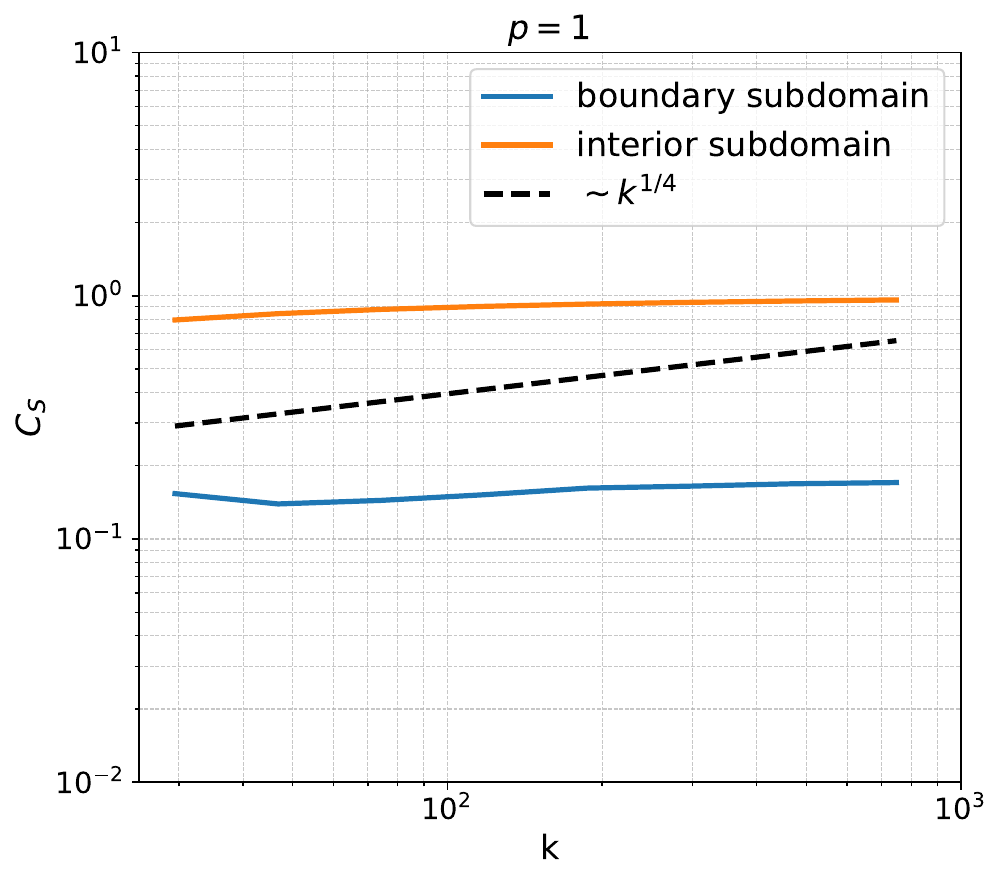}\qquad
    \includegraphics[scale=0.43]{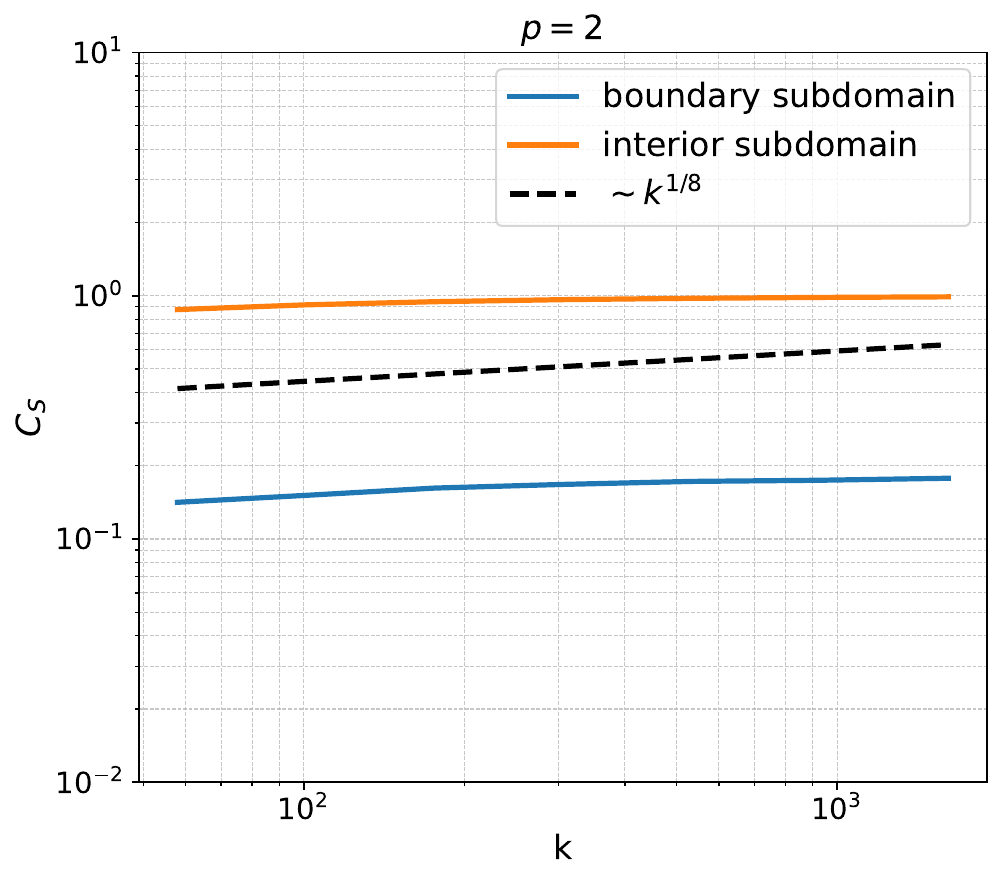}
    \caption{The quantity $C_S$ for an interior and a boundary subdomains for linear (left) and quadratic (right) finite elements.}
    \label{fig:stability}
\end{figure}

Next, we investigate the dependence of the stability estimate \eqref{eq:error_local_Helmholtz} on the term $(kh)^{-1/2}$. To do the test, we set the right hand side $f=1$, and choose the boundary data $g=\nabla \cm{u^{\mathdutchcal{e}}}\cdot {\bm n} - {\rm i}k \cm{u^{\mathdutchcal{e}}}$ with $\cm{u^{\mathdutchcal{e}}} = \exp\big({\rm i}\bm{k} \cdot {\bm x}\big)$ as above. In \cref{fig:stability}, we fix $kH_i^{\ast}$ and compute the quantity 
\begin{equation*}
C_{S} = 
\frac{\Vert u^{{\mathdutchcal{e}}}_{h}- \psi_{h,i}\Vert_{\mathcal{A},\omega_{i}^{\ast},k}} 
{\Vert u^{{\mathdutchcal{e}}}_{h}\Vert_{\mathcal{A},\widetilde{\omega_{i}^{\ast}},k}}
\end{equation*}
for an interior and a boundary subdomains for linear and quadratic finite elements. With this setup, estimate \eqref{eq:error_local_Helmholtz} gives the theoretical bound $C_{S} \lesssim  (kh)^{-1/2} \sim k^{\frac{1}{4p}}$. However, as can be seen from \cref{fig:stability}, this theoretical bound is rather pessimistic.



\begin{figure}[h]
    \centering
    \includegraphics[scale=0.40]{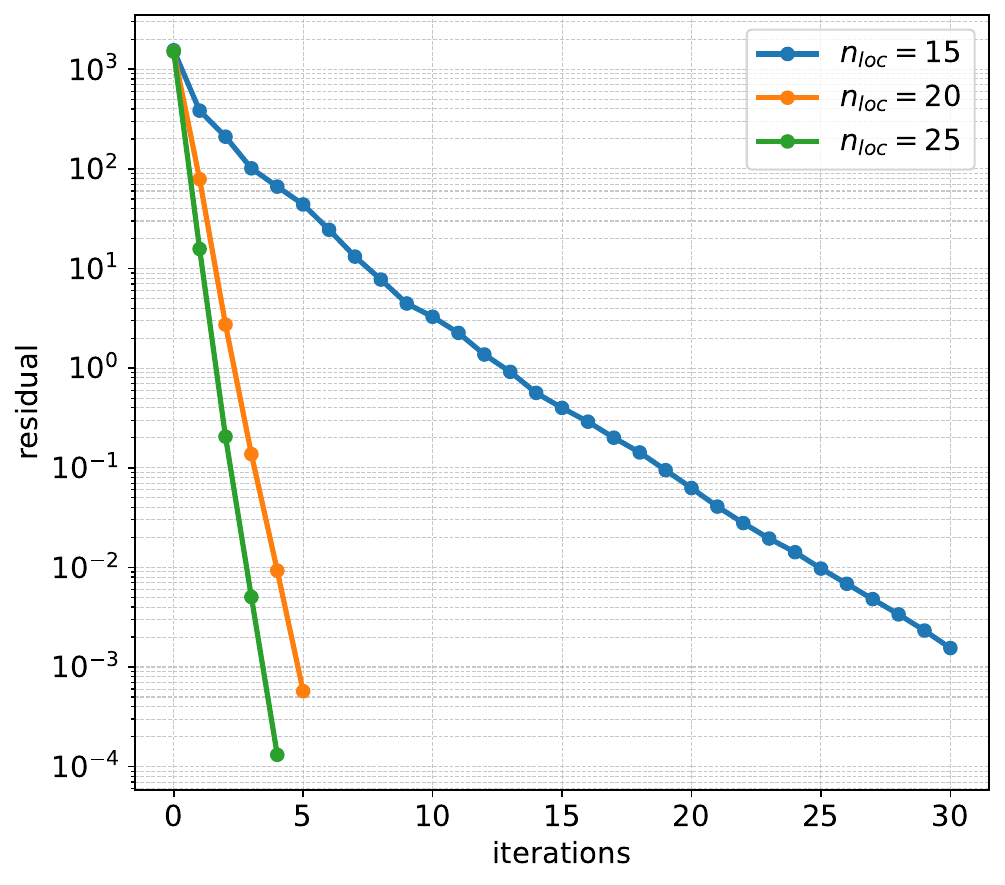}
     \includegraphics[scale=0.40]{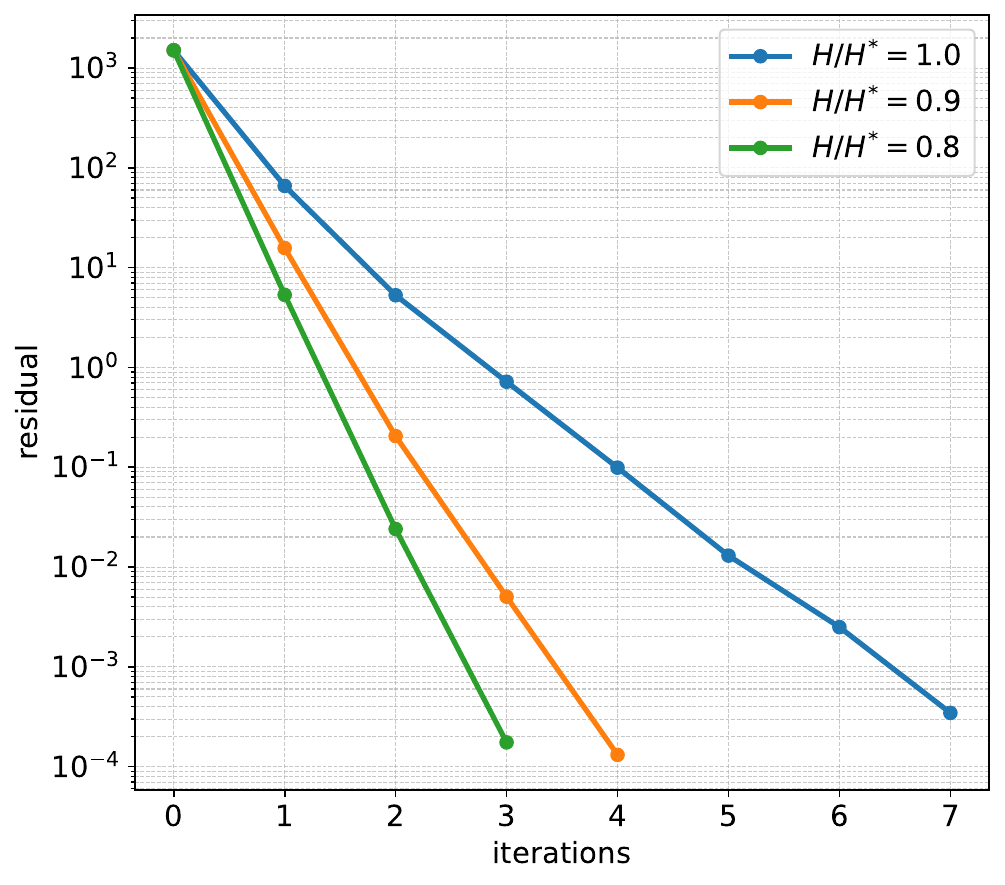}
    \caption{Constant-coefficient example \RS{with $k=200$, $h=1/752$ and $m=16$}: Residuals vs. iteration numbers for variable $n_{loc}$ (left, with $H/H^{\ast}=0.9$) and oversampling ratios (right, with $n_{loc}=25$).}
    \label{fig:residuals}
\end{figure}




\begin{table}[h]
\centering
\begin{tabular}{lll|ll|ll}
\toprule
 & \multirow{2}{*}{$k$} & \multirow{2}{*}{\RS{$M$}} &\multicolumn{2}{c}{$H/H^*=1.0$} & \multicolumn{2}{c}{$H/H^*=0.9$} \\
 & & &$\rho =2 k^{-0.5}$ & $\rho = 20 k^{-1}$ & $\rho =2  k^{-0.5}$ & $\rho =20 k^{-1}$ \\
\midrule
\multirow{4}{*}{$H\sim k^{-0.4}$} &200 &17$\times$17 & 4 (9567) & 3 (11968) & 4 (5440) & 3 (6239) \\
 &400 &22$\times$22 & 3 (23624) & 3 (38324) & 4 (15184) & 3 (18880) \\
 &600 &26$\times$26 & 3 (40924) & 3 (75976) & 4 (23976) & 3 (29176) \\
 &800 &29$\times$29 & 3 (58464) & 2 (122583) & 4 (34779) & 3 (41271) \\
 \midrule
 \multirow{4}{*}{$H\sim k^{-0.6}$}&200 &19$\times$19 & 4 (10655) & 3 (13680) & 4 (6479) & 3 (7129) \\
 &400 &29$\times$29 & 4 (31751) & 3 (51968) & 4 (20325) & 3 (25143) \\
 &600 &37$\times$37 & 4 (58608) & 3 (110519) & 4 (35927) & 3 (43849) \\
 &800 &44$\times$44 & 4 (88884) & 3 (187096) & 5 (54912) & 3 (73572) \\
\midrule
\multirow{4}{*}{$H\sim k^{-1.0}$} &200 &25$\times$25 & 4 (15025) & 4 (19200) & 4 (8975) & 4 (9600) \\
&400 &50$\times$50 & 4 (56304) & 3 (88200) & 5 (36700) & 3 (46500) \\
&600 &76$\times$76 & 5 (119776) & 3 (233476) & 6 (85424) & 3 (119476) \\
&800 &101$\times$101 & 6 (212201) & 3 (474199) & 10 (151799) & 3 (232199) \\
\bottomrule
\end{tabular}
\caption{Constant-coefficient example: \#GMRES iterations for different choices of subdomain size, eigenvalue tolerance and oversampling. The size of \RS{the} coarse space is shown in parentheses. The average sizes of the fine-scale FE problems are: $1.5\times 10^{5}$, $8.8\times 10^{5}$, $2.4\times 10^{6}$, $4.9\times 10^{6}$.}
\label{tab:robin_tau}
\end{table}

\begin{table}[h]
\parbox{.45\linewidth}{
\centering
\begin{tabular}{lrrr}
\toprule
$k\textbackslash H/H^*$ & 1.0 & 0.9 & 0.8 \\
\midrule
200 & 10 & 8 & 10 \\
400 & 8 & 5 & 5 \\
600 & 6 & 4 & 4 \\
800 & 6 & 4 & 3 \\
1000 & 6 & 4 & 3 \\
\bottomrule
\end{tabular}

\label{tab:iterations_increasing_wavenumber_and_ndom_2}
}
\parbox{.45\linewidth}{
\centering
\begin{tabular}{lrrr}
\toprule
$k\textbackslash H/H^*$ & 1.0 & 0.9 & 0.8 \\
\midrule
200 & 7 & 5 & 4 \\
400 & 6 & 4 & 3 \\
600 & 6 & 3 & 2 \\
800 & 5 & 3 & 2 \\
1000 & 6 & 3 & 2 \\
\bottomrule
\end{tabular}

}
\caption{Constant-coefficient example: \#GMRES iterations for case \RS{(i)}  with $m=k/25$, $n_{loc}=4\log_2(k)$, and $\delta = 1/(10m)$ \RS{(left)}, and for case \RS{(ii)} with $m=16$, $n_{loc} = k/10$, and $\delta = 2h$ \RS{(right)}.}
\label{tab:iterations_increasing_wavenumber_7}
\end{table}

\begin{table}[h]
\parbox{.45\linewidth}{
\centering
\begin{tabular}{lrrr}
\toprule
$h\textbackslash n_{loc}$ & 15 & 20 & 25 \\
\midrule
1/200 & 12 & 4 & 3 \\
1/400 & 12 & 4 & 3 \\
1/600 & 12 & 4 & 3 \\
1/800 & 12 & 4 & 3 \\
\bottomrule
\end{tabular}
}
\parbox{.45\linewidth}{
\hspace{5ex}
\centering
\begin{tabular}{lrrr}
\toprule
$ \RS{M} \textbackslash n_{loc}$ & 10 & 15 & 20 \\
\midrule
8$\times$8 & 65 & 13 & 5 \\
16$\times$16 & 9 & 4 & 3 \\
32$\times$32 & 4 & 3 & 3 \\
64$\times$64 & 3 & 3 & 2 \\
\bottomrule
\end{tabular}
}
\caption{Constant-coefficient example: \#GMRES iterations for decreasing mesh-size (left) with $m = 8$, $k = 100$, $\delta=\delta^{\ast} = 1/100$, and for increasing number of subdomains (right) with $h=1/512$, $k = 100$, $\delta=\delta^{\ast} = 1/(8m)$.}
\label{tab:robustness_meshsize_20241207_124144}
\end{table}




\subsubsection{Performance of MS-GFEM Preconditioner}
\label{subsubsect:constant_example_preconditioner}

In this subsection, we evaluate the performance of MS-GFEM as a preconditioner within GMRES. For all tests, we choose the data $g$ and $f$ such that the problem \eqref{eq:1-1} has the exact solution $u^{\mathdutchcal{e}} = \exp\big({\rm i}\bm{k} \cdot {\bm x}\big)$ with $\bm{k} = k/\sqrt{2}(1,\, 1)$, and we use quadratic elements with $h\sim k^{-\frac54}$. In \cref{fig:residuals}, we fix the wavenumber and compare the convergence performance of GMRES for different $n_{loc}$ (the number of eigenfunctions used per subdomain) and oversampling ratios. The results clearly show that using more eigenfunctions or oversampling leads to significantly faster convergence. Remarkably, the residual reduction that the method achieves with $n_{loc} = 20$ in 5 iterations is smaller than what it achieves with $n_{loc} = 15$ in 30 iterations.

Next, we investigate the iteration counts and size of the coarse 
space for increasing $k$. Here we use an eigenvalue tolerance $\rho$ to select eigenfunctions on each subdomain -- all eigenfunctions corresponding to eigenvalues that are larger than $\rho^2$ are used \RS{($\lambda_{\mathrm{max}}^{1/2}\leq \rho$)}. Moreover, we scale the size $H$ of the subdomains according to $H\sim k^{-\tau}$ with $\tau>0$. A larger $\tau$ leads to a smaller subdomain size and a larger number of subdomains. \Cref{tab:robin_tau} gives the numbers of GMRES iterations and sizes of coarse spaces for different numbers of $\tau$ and different eigenvalue tolerances $\rho$. We see that for a given choice of $\tau$, $\rho$, and oversampling ratio $H/H^{\ast}$, the iteration counts remain roughly the same \RS{when} increasing $k$. However, the sizes of coarse spaces are noticeably different. \RS{A larger $\tau$ leads to a markedly faster growth in the coarse problem size as $k$ increases despite a smaller local problem size. The choice of $\rho = 2k^{-0.5}$ is more efficient than $\rho = 20k^{-1}$ -- it results in comparable iteration counts with much smaller coarse spaces, especially when no oversampling is used. Using oversampling significantly reduces the coarse space size, especially for the choice of $\rho = 20k^{-1}$. It is worth noting that while only the parameter setting ($H\sim k^{-1}$, $\rho \sim k^{-1}$) is covered by our theory, the method with the other settings also performs very well.}



We have learned from the frequency scaling tests, how 
the number of eigenfunctions needs to be chosen depending on $k$ to achieve a given 
accuracy. Now, we use these scalings of $n_{loc}$ in the iterative method and display the iterations counts in \cref{tab:iterations_increasing_wavenumber_7} for two choices of $H$ and $n_{loc}$: \RS{(i)} $H\sim k^{-1}$, $n_{loc}\sim \log_{2}(k)$ and \RS{(ii)} $H$ fixed, $n_{loc}\sim k$. We see that for \RS{both} cases, the iteration counts are comparable and remain roughly the same or even decrease as $k$ increases. Moreover, using oversampling reduces the number of iterations required. 

Finally, we test the robustness of the preconditioner for decreasing fine mesh-size $h$ and increasing number of subdomains $m$. The results are shown in \cref{tab:robustness_meshsize_20241207_124144}. We see that the iteration counts remain unchanged as we decrease the fine mesh-size. On the other hand, the method converges faster as we increase the number of subdomains with fixed $n_{loc}$. This can be explained as follows: as the size of subdomains becomes smaller, there are fewer `waves' in each subdomain, and the solution can be better and better resolved by a fixed number of local eigenfunctions. 

\subsection{Marmousi Model}\label{subsec:Marmousi model}
\begin{figure}
\centering
\includegraphics[scale=0.28]{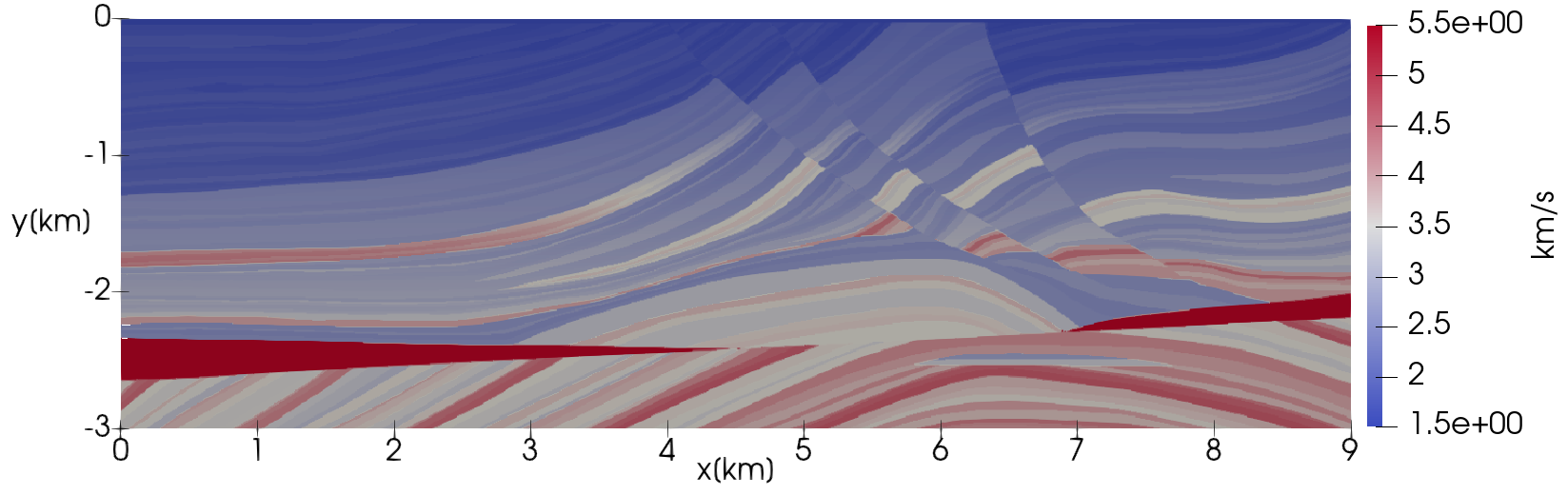}
\caption{The velocity field of the Marmousi model.}
\label{fig:6-6}
\end{figure}

\RS{We now} test our method for \RS{heterogeneous problems, starting with} the Marmousi model \cite{versteeg1994marmousi}. The model is a widely used benchmark in seismic imaging and geophysics, specifically designed to simulate complex subsurface structures. We use the same problem setup as in \cite[Section 6.3]{ma2023wavenumber} which we recap here. The problem is posed on the domain $(0, \,9\,{\rm km})\times (-3\,{\rm km}, 0)$ and a point source is placed near the top boundary. Homogeneous Dirichlet and impedance boundary conditions are imposed on the top side and on the remaining three sides, respectively. Moreover, $A({\bm x})=I$ \RS{and} $\beta({\bm x})=kV({\bm x})$ \RS{with} the velocity field $1/V({\bm x})$ is shown in \cref{fig:6-6}. The computational setting is the same as in the preceding subsection, except that we choose $m$ and $3m$ subdomains in the $y$ and $x$ directions, respectively. 

In \cref{tab:hgeneo_table20241211_164637}, we display the iteration counts for increasing frequencies $\nu$, \RS{varying the} oversampling ratios $H/H^{\ast}$, \RS{the} numbers of eigenfunctions used per subdomain $n_{loc}$, and \RS{the} numbers of subdomains \RS{$M=3m^2$}. \RS{The symbol $\times$ represents that the method did not converge in 100 iterations.} For all the tests, we use quadratic finite elements with a fixed mesh-size $h=7.5m$, corresponding to a resolution with 10 points per (minimal) wavelength for all the frequencies. We observe that for $\nu = 5$\,Hz and $\nu = 10$\,Hz, the method performs very well in terms of iteration numbers for all the choices of \RS{$M$} and $n_{loc}$, and that the increase in the iteration number from $\nu = 5$\,Hz to $\nu = 10$\,Hz is mild. However, the iteration numbers significantly increase for the higher frequency $\nu = 20$\,Hz when relatively small $n_{loc}$ and \RS{$M$} are used. But for larger $n_{loc}$ \RS{or $M$}, the iteration numbers for $\nu = 20$\,Hz remain reasonably small. This observation indicates that for a given frequency, there is a threshold number of eigenfunctions needed for the method to converge, and that decreasing subdomain size can reduce the threshold, which agrees with our observation in the frequency scaling tests. In general, we observe a reduced iteration count if 
oversampling is used. 


\begin{table}[h!]
\centering
\begin{tabular}{llllll|llll}
\toprule
 &  & \multicolumn{4}{c}{$\cB H/H^{\ast}=1$} & \multicolumn{4}{c}{$\cB H/H^{\ast}=0.86$} \\[0.5ex]
 $\cB \nu $ & $\cB \RS{M} \textbackslash n_{loc}$ & 15 & 20 & 25 & 30 & 15 & 20 & 25 & 30 \\
\midrule
 & 300  & \cR 5 & \cR 4 &\cR 4 &\cR 3 &\cR 4 &\cR 3 &\cR 2 &\cR 2\\
\multirow[t]{3}{*}{5\,Hz} & 1200 &\cR 5 &\cR 4 &\cR 3 &\cR 3 &\cR 3 &\cR 2 &\cR 2 &\cR 2\\
 & 4800 &\cR 4 &\cR 3 &\cR 3 &\cR 3 &\cR 3 &\cR 2 &\cR 2 &\cR 2 \\
 \cline{1-10}
  & 300 &\cR 12 &\cR 7 &\cR 5 &\cR 4 &\cR 4 &\cR 3 & \cR 2 &\cR 2\\
\multirow[t]{3}{*}{10\,Hz} &1200 &\cR 7 &\cR 5 &\cR 4 &\cR 4 &\cR 4 &\cR 3 &\cR 2 &\cR 2\\
 & 4800 &\cR 5 &\cR 4 &\cR 4 &\cR 3 &\cR 3 &\cR 2 &\cR 2 &\cR 2 \\
 \cline{1-10}
   &300&\cR $\times$ &\cR $\times$ &\cR 30 &\cR 7 &\cR $\times $ & \cR$\times $ & \cR 24 &\cR 5 \\
\multirow[t]{3}{*}{20\,Hz} & 1200 &\cR 28 &\cR 11 &\cR 6 &\cR 5 &\cR 22 &\cR 5 &\cR 4 &\cR 3 \\
 & 4800 &\cR 10 &\cR 6 &\cR 4 &\cR 4 &\cR 6 &\cR 4 &\cR 3 &\cR 3 \\
\bottomrule
\end{tabular}

 \caption{Marmousi model: \#GMRES iterations for different 
frequencies $\nu$, oversampling ratios $H/H^{\ast}$, numbers of eigenfunctions used per subdomain $n_{loc}$, and numbers of subdomains in the $y$-direction $m$, with $\delta=40h/m$.}
 \label{tab:hgeneo_table20241211_164637}
\end{table}



\subsection{High Contrast Example}

In this example, we assess the performance of the preconditioner for high-contrast coefficients. The problem setup is the same as in section 6.1 of \cite{peterseim2020computational} (see also \cite{verfurth2024numerical}): $\Omega=(0,1)^2$, $V({\bm x}) = 1$, $k=9$, $\beta({\bm x}) =1$, $g({\bm x}) = 0$, and $f$ is a localized source at $(0.125, 0.5)$. For $\varepsilon >0$, the piecewise-constant coefficient $A$ is 
given by 
\begin{equation*}
    A({\bm x}) = \begin{cases}
    \varepsilon^2, \quad \textit{ if } \bm{x} \in \Omega_{\varepsilon},\\
    1, \quad \;\;\textit{ else},
    \end{cases}
\end{equation*}
where 
\begin{equation*}
    \Omega_\varepsilon=(0.25, 0.75)^2\cap \bigcup_{j\in \mathbb{Z}^2}\varepsilon(j+(0.25, 0.75)^2),
\end{equation*}
namely, $\Omega_{\varepsilon}$ is the union of periodic inclusions located in the center of $\Omega$; see \cref{fig:periodic_inclusions} for an illustration of $\Omega_{\varepsilon}$ with $\varepsilon = 1/4$. The problem setup is of physical interest, \RS{as it} can result in resonance behavior due to the presence of \RS{the} small inclusions. We refer to \cite{peterseim2020computational,verfurth2024numerical} for a detailed discussion. For $0<\varepsilon \ll 1$, the problem is considered hard, because (i) the coefficient $A$ is strongly heterogeneous with a high contrast and (ii) the effective wavenumber $k/\varepsilon$ in the inclusions is large. 

We discretize the problem with quadratic finite elements, and set the mesh-size $h\approx 2(k/\varepsilon)^{-5/4}$ to ensure a reasonable resolution for the solution in the inclusions. The rest of the computational setting is the same as in \cref{subsubsect:constant_example_preconditioner}. In particular, we use a prescribed eigenvalue tolerance $\rho$ for selecting the eigenfunctions as in \cref{subsubsect:constant_example_preconditioner}.  \Cref{tab:high_contrast_epsilon_scaling20241220_170601} shows the 
iteration counts and sizes of the coarse spaces for decreasing $\varepsilon$ and various eigenvalue tolerances $\rho$ with fixed subdomains. We see that for a wide range of \RS{tolerances} $\rho$, the iteration counts stay bounded or even decrease with decreasing $\varepsilon$. This demonstrates that for this problem, our coarse space can automatically select physically relevant local modes such that the resulting preconditioner is robust to the parameter $\varepsilon$. \RS{Moreover}, we observe that for a fixed tolerance $\rho$, the size of the coarse space grows \RS{only slightly faster than} linearly with $\varepsilon^{-1}$ -- a moderate growth compared with the fine-scale problem whose size grows as $O(\varepsilon^{-\frac52})$.  

\begin{figure}
    \centering
    \begin{tikzpicture}
        \draw[thick] (0,0) rectangle (4,4);
        \node[left] at (0.8,0.5) {\Large $\Omega$};


        \foreach \x in {1.125,1.625,2.125,2.625}
            \foreach \y in {1.125,1.625,2.125,2.625}
                \filldraw[fill=red!75] (\x,\y) rectangle ++(0.25,0.25);

        \draw[thick] (2.5,2.5) rectangle ++(0.5,0.5);

        \draw[<->] (2.5,3.1) -- (3.0,3.1);
        \node[above] at (2.75,3.1) {\small $\varepsilon$};

        \draw[thick] (4.5,1) rectangle (6.5,3);
        \node[above] at (5.5,3) {\Large $Y$};

        \filldraw[fill=red!75] (5,1.5) rectangle (6,2.5);

        \draw[thick] (3.0,2.8) -- (4.5,3);
        \draw[thick] (2.9,2.5) -- (4.5,1.5);
    \end{tikzpicture}
\caption{Illustration of $\Omega_{\varepsilon}$ with $\varepsilon = \frac{1}{4}$. The rescaled unit cell is shown on the right.}    
\label{fig:periodic_inclusions}
\end{figure}

\begin{table}
\centering
\begin{tabular}{lllll}
\toprule
$\rho\textbackslash \varepsilon $& $2^{-3}$ & $2^{-4}$ & $2^{-5}$ & $2^{-6}$ \\
\midrule
$2^{-1}$ & 52 (416) & 15 (1248) & 21 (3168) & 22 (7112) \\
$2^{-2}$ & 15 (1108) & 6 (2036) & 10 (4904) & 9 (11952) \\
$2^{-3}$ & 9 (1536) & 5 (3420) & 5 (7252) & 3 (17380) \\
$2^{-4}$ & 8 (2076) & 4 (4768) & 4 (10140) & 2 (23648) \\
$2^{-5}$ & 5 (3200) & 3 (6288) & 2 (13252) & 2 (31312) \\
\bottomrule
\end{tabular}

\caption{High contrast example: \#GMRES iterations and size of coarse space (in parentheses) with $m = 20$, $\delta = \delta^{\ast}=2h$. $\rho$ denotes the eigenvalue tolerance.}
\label{tab:high_contrast_epsilon_scaling20241220_170601}
\end{table}

\subsection{3D EAGE/SEG Overthrust Model}
\begin{figure}
    \centering
    \includegraphics[scale =0.33]{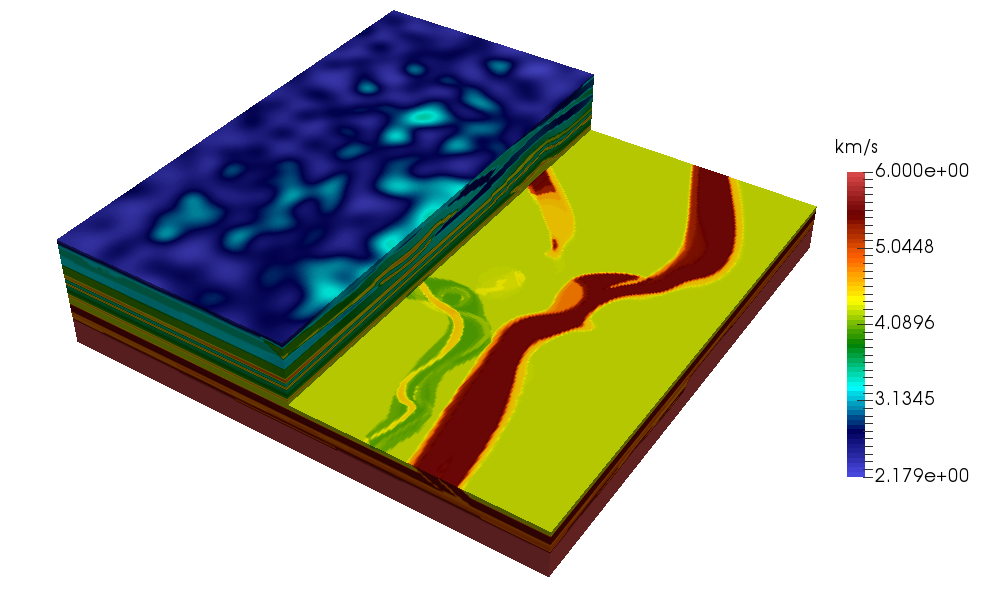}\hspace{4ex}
     \includegraphics[scale =0.65]{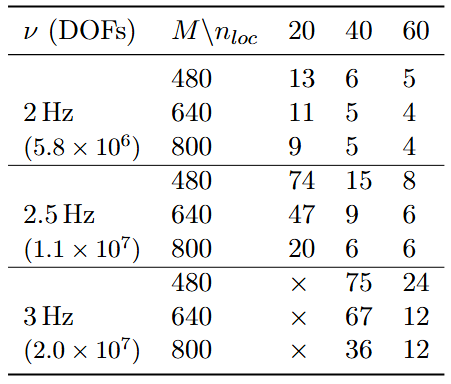}
    \caption{3D Overthrust model: the velocity field of the model (left) and \#GMRES iterations (right) for different $M$ (number of subdomains) and $n_{loc}$ (number of eigenfunctions per subdomain). The size of the fine-scale problem is given in the parentheses.  }
    \label{fig:overthrust_example}
\end{figure}

In this example, we test our method on the 3D SEG/EAGE Overthrust model \cite{aminzadeh19963} to evaluate its capability for three-dimensional large-scale problems. The model can be seen as the three-dimensional counterpart of the Marmousi model, spanning 20 × 20 × 4.65 km with a heterogeneous velocity field as depicted in \cref{fig:overthrust_example} (left). We use the same problem setup as in \cite[section 6.3]{bootland2021comparison}, where homogeneous Dirichlet and impedance boundary conditions are imposed on the top face and the remaining five faces, respectively, and a point source is placed at (2.5, 2.5, 0.58) km. We perform simulations using P2 finite elements on tetrahedral meshes for frequencies of 2\,Hz, 2.5\,Hz, and 3\,Hz. For each frequency, we use 8 points per (minimal) wavelength for the discretization, which is a reasonable resolution in practice. For domain partition, we use the automatic graph partitioner \textit{Metis} \cite{karypis1997metis}. The overlapping subdomains and oversampling domains are constructed by adding one and two layers of adjoining mesh elements to the non-overlapping subdomains, respectively, \RS{such that} $\delta = \delta^{\ast} = h$. The construction of partition of unity functions follows \cite{jolivet2014overlapping}. In \cref{fig:overthrust_example} (right), we display the iteration counts for different numbers of subdomains $M$ and numbers of eigenfunctions used per subdomain $n_{loc}$. \RS{The maximum number of iterations was set to 100}. We observe that the method converges faster with increasing $M$ and $n_{loc}$. As for the Marmousi model, the method fails to converge for a high frequency if $n_{loc}$ is too small, but it performs well if a moderate number of eigenfunctions are used. Notably, in the case of $\nu = 3\,$Hz, the method is able to converge in 12 iterations with a fairly small coarse problem (size 38400) -- using 640 subdomains and 60 eigenfunctions per subdomain. By contrast, the size of the fine-scale problem is about $2.08 \times 10^{7}$.

\section{Conclusions}\label{sec:conclusion}
We have formulated the MS-GFEM for Helmholtz problems as a two-level hybrid, restricted additive Schwarz preconditioner. Based on the theory of MS-GFEM developed in our previous work, we give a detailed characterization of the convergence rate of the preconditioned GMRES algorithm. Numerical results show that the preconditioner can achieve fast convergence with a relatively small coarse space, by flexibly adjusting the subdomain size and/or the number of local eigenfunctions used. This desirable performance makes the method \RS{particularly} appealing for seismic imaging applications with \RS{highly varying wave speeds and} a large number of sources.

There are two practical issues concerning the preconditioner that need further investigation. The first is the solution of the local eigenproblems, which makes up the bulk of the computational work in a single simulation. While the implicitly restarted Arnoldi algorithm implemented in Arpack is efficient for the reformulated eigenproblems \eqref{mixed_EVP} in terms of computational time, its applicability is limited by available memory, which consequently limits the subdomain size. Alternative eigensolvers that have a lower memory requirement include the Jacobi-Davidson method \cite{sleijpen2000jacobi} and the LOBPCG method \cite{knyazev2001toward}, but they may need a good preconditioner to speed up convergence. The selection or design of an efficient eigensolver for our method is an important practical issue and the investigation is now ongoing. The second issue is the solution of the coarse problem. Indeed, the coarse problem needed in our preconditioner is relatively small in size, which can be handled efficiently by direct solvers for fine-scale problems with up to tens of millions of unknowns. However, for larger-scale problems, the factorization of the coarse problem matrix may be costly and an iterative coarse solver will inevitably be needed. In the near future, we will investigate extending the proposed two-level preconditioner to a multilevel Schwarz preconditioner.

\vspace{2ex}
\textbf{Acknowledgements.} The authors thank Ivan G. Graham (University of Bath), Euan A. Spence (University of Bath), Victorita Dolean (Eindhoven University of Technology), Shihua Gong (The Chinese University of Hong Kong, Shenzhen), and Jens Markus Melenk (Vienna University of Technology) for many stimulating discussions. In particular, the authors thank IGG for suggesting the frequency scaling tests. 
\bibliographystyle{siamplain}
\bibliography{references}

\end{document}